\documentstyle[amssymb,12pt]{amsart}
\newtheorem{prop}{Proposition}[section]
\newtheorem{prop:def}{Proposition-Definition}[section]

\newtheorem{lemma}{Lemma}[section]
\newtheorem{thm}{Theorem}[section]
\newtheorem{cor}{Corollary}[section]

\theoremstyle{remark}
\newtheorem{remark}{Remark}

\begin{document}
\newcommand{\nc}{\newcommand} \nc{\on}{\operatorname}
\nc{\pa}{\partial} \nc{\cA}{{\cal A}}\nc{\cB}{{\cal B}}\nc{\cC}{{\cal
    C}} \nc{\cE}{{\cal E}}\nc{\cG}{{\cal G}}\nc{\cH}{{\cal H}}
\nc{\cX}{{\cal X}}\nc{\cR}{{\cal R}}\nc{\cL}{{\cal L}} \nc{\cK}{{\cal
    K}} \nc{\sh}{\on{sh}}\nc{\Id}{\on{Id}}\nc{\Diff}{\on{Diff}}
\nc{\ad}{\on{ad}}\nc{\Der}{\on{Der}}\nc{\End}{\on{End}}
\nc{\res}{\on{res}}\nc{\ddiv}{\on{div}}
\nc{\card}{\on{card}}\nc{\dimm}{\on{dim}}
\nc{\Jac}{\on{Jac}}\nc{\Ker}{\on{Ker}}
\nc{\Imm}{\on{Im}}\nc{\limm}{\on{lim}}\nc{\Ad}{\on{Ad}}
\nc{\ev}{\on{ev}} \nc{\Hol}{\on{Hol}}\nc{\Det}{\on{Det}}
\nc{\Bun}{\on{Bun}}\nc{\diag}{\on{diag}}
\nc{\de}{\delta}\nc{\si}{\sigma}\nc{\ve}{\varepsilon}\nc{\z}{\zeta}
\nc{\al}{\alpha}\nc{\vp}{\varphi} \nc{\CC}{{\mathbb
    C}}\nc{\ZZ}{{\mathbb Z}} \nc{\NN}{{\mathbb N}}\nc{\zz}{{\mathbf
    z}} \nc{\AAA}{{\mathbb A}}\nc{\cO}{{\cal O}} \nc{\cF}{{\cal
    F}}\nc{\cM}{{\cal M}} \nc{\la}{{\lambda}}\nc{\G}{{\mathfrak
    g}}\nc{\mm}{{\mathfrak m}} \nc{\A}{{\mathfrak a}}
\nc{\HH}{{\mathfrak h}} \nc{\N}{{\mathfrak n}}\nc{\B}{{\mathfrak b}}
\nc{\La}{\Lambda}
\nc{\g}{\gamma}\nc{\eps}{\epsilon}\nc{\wt}{\widetilde}
\nc{\wh}{\widehat} \nc{\bn}{\begin{equation}}\nc{\en}{\end{equation}}
\nc{\SL}{{\mathfrak{sl}}}\nc{\ttt}{{\mathfrak{t}}}

%
%
%

\newcommand{\ldar}[1]{\begin{picture}(10,50)(-5,-25)
\put(0,25){\vector(0,-1){50}}
\put(5,0){\mbox{$#1$}} 
\end{picture}}

\newcommand{\lrar}[1]{\begin{picture}(50,10)(-25,-5)
\put(-25,0){\vector(1,0){50}}
\put(0,5){\makebox(0,0)[b]{\mbox{$#1$}}}
\end{picture}}

\newcommand{\luar}[1]{\begin{picture}(10,50)(-5,-25)
\put(0,-25){\vector(0,1){50}}
\put(5,0){\mbox{$#1$}}
\end{picture}}

\title[Commuting differential operators associated to complex curves]
{Commuting differential and difference operators associated to complex
  curves I}

\author{B. Enriquez}

\address{B.E.: Centre de Math\'ematiques, Ecole Polytechnique, 
UMR 7640 du CNRS, 91128 Palaiseau, France}

\address{FIM, ETH-Zentrum, HG G46, CH-8092 Z\"urich, Switzerland}

\author{G. Felder}

\address{G.F.: D-Math, ETH-Zentrum, HG G44, CH-8092 Z\"urich,
  Switzerland}

\date{July 1998}

\maketitle

\subsection*{Introduction}
In \cite{FS}, B. Feigin and A. Stoyanovsky introduced a remarkable
parametrization of the space of conformal blocks associated with a
complex curve $X$ and a semisimple Lie algebra $\bar\G$. This space is
defined as the set of $\G^{out}$-invariant forms on an integrable
$\G$-module $L_{\La,k}$ located at a point $P_0$ of $X$, $\G^{out}$
being the Lie algebra of regular maps from $X - \{P_0\}$ to $\bar\G$
and $\G$ the Kac-Moody Lie algebra at $P_0$. Feigin and Stoyanovsky
associate to such a form $\psi$, the family of forms on a product of
symmetric products of $X$
\begin{equation} \label{FeSt} 
  f_{FS}(z^{(i)}_j) = \langle \psi, \prod_{i=1}^r e_i(z^{(i)}_j)
  dz^{(i)}_j (wv_{top}^{(P_0)}) \rangle ,
\end{equation} 
where $r$ is the rank of $\bar\G$, $e_i(z)dz$ are the currents
associated with the simple nilpotent generators of $\bar\G$, $w$ is an
affine Weyl group element and $v_{top}^{(P_0)}$ is the highest weight
vector of $L_{\La,k}$.

In this paper, we introduce the twisted conformal blocks $\psi_\la =
\psi\circ e^{\sum_{i,a} \la^{(i)}_a h_i[r_a]}$. Here $(r_a)_{a = 1,
  \ldots,g}$ are functions on $X$, regular outside $P_0$,
single-valued around $a$-cycles and all $b$-cycles except the $a$th,
along which it increases by $1$ (sect. \ref{bases}), and the $h_i$ are
the simple coroots of $\bar\G$. The functions $r_a$ are thus the
analogues of the function $\theta'/\theta$ in the elliptic case.
$\psi_\la$ is independent of the choice of the functions $r_a$.

For any $v$ in $L_{\La,k}$, the function $v\mapsto \langle
\psi_\la,v\rangle$ is defined as a formal function in $\la =
(\la^{(i)}_a)$. We show (Thm.\  \ref{seerose}) that it is actually a
holomorphic function in $\la$ with theta-like properties.  This result
relies on adelization of the representations $L_{\La,k}$ (see
\cite{TUY}), reduction to the $\SL_2$ case, formulas for the tame
symbol and the identity $(f^k)' = - :hf^k:$ (see \cite{Lep-Primc}).
This generalizes a result obtained in \cite{FW} in the genus $1$ case.

We then consider the forms
\begin{equation} \label{oo}
f_\la(z^{(i)}_j) = \langle \psi_\la, \prod_{i=1}^r e_i(z^{(i)}_j)
dz^{(i)}_j (wv_{top}^{(P_0)}) \rangle. 
\end{equation} 
These forms have the following geometric interpretation.  It is known
(\cite{BL,Kumar}) that conformal blocks can be viewed as sections of a
bundle on the moduli space $Bun_{\bar G}$ of $\bar G$; such sections
are called generalized theta-functions.  In sect.  \ref{lifts}, we
explain that the forms (\ref{FeSt}) of Feigin-Stoyanovsky can be
viewed as generating functions for lifts of the generalized
theta-functions to a space, which in the case $\bar\G = \SL_n$ can be
described as $Bun_{(n_i,P_0)}$, the moduli space of bundles with
filtration $E_1\subset E_2 \subset \cdots $ and associated graded
isomorphic to $\oplus_i \cO(n_i P_0)$, $n_i$ some integer numbers.
From this viewpoint, the twisted correlation functions (\ref{oo}) are
generating functions for lifts of generalized theta-functions to the
moduli space $Bun_B$ of $B$-bundles over $X$, where $B$ is the Borel
subgroup of $\bar G$.

We then express the Knizhnik-Zamolodchikov-Bernard (KZB) connection in
terms of the forms (\ref{oo}) (sect.  \ref{KZB}). Our treatment of the
KZB connection follows \cite{Houches}; the KZB connection is defined
on the space of projective structures on curves of genus $g$.
However, we such a projective structure is canonically attached to the
choice of $a$-cycles on the curve, via a bidifferential form
$\wt\omega$ (see (\ref{wt:omega}); this form appeared in \cite{Fay},
cor.\ 2.6). This allows to define the KZB connection as a projectively
flat connection on the moduli space of curves with marked $a$-cycles,
which is intermediate between the moduli space of curves and its
universal cover.  The KZB connection is expressed as the action of
differential-evaluation operators $(T_z)_{z\in X}$ on the
$f_\la(z_j^{(i)})$, which are forms on $J^0(X)^r \times \prod_i
S^{n_i}X$ (differential in $\la$ and residues and evaluation in the
$\la^{(i)}_a$).

We also express the KZB connection in the directions given by
variation of points in case of a curve with marked points (sect.
\ref{marked}). Denote by $\wt m$ a quadruple
$(X,[\{\zeta_i\}],P_i,\zeta_i)$ formed by a curve with projective
structure, marked points and coordinates at these points, by $\psi(\wt
m)$ a conformal block associated to this complex structure, and by
$f(\wt m)_\la(z_\al)$ the twisted correlation function associated with
this conformal block according to (\ref{oo}).
In the case $\bar\G = \SL_2$, the connection takes the form
$$ 
  2(k+2)\nabla_{{{\pa}\over{\pa P_i}}} f(\wt m)_\la(z_\al) = 2(k+2)
  {{\pa}\over{\pa P_i}} f(\wt m)_\la(z_\al) - (K_if)(\wt m)_\la(z_\al), 
$$
with 
\begin{align} \label{cite:KZ} 
  & (K_if)(\wt m)_\la(z_\al) = [ - \La_i\sum_a
  \omega_a(P_i)\pa_{\la_a} + \La_i \left( \sum_{j\neq i} \La_j
    G(P_j,P_i) - 2 \sum_\al G(z_\al,P_i) \right) \\ & \nonumber +
  \La_i^2 \phi(P_i) + 2\La_i g_{2\la}(P_i) ] f(\wt m)_\la(z_\al) \\ &
  \nonumber + \sum_{\al}[-2 G_{2\la}(P_i,z_\al)(\sum_a \omega_a(z_\al)
  \pa_{\la_a} + 2 \sum_{\beta \neq \al} G(z_{\al},z_\beta)) \\ &
  \nonumber -4 G_{2\la}(P_i,z_\al)G(z_\al,P_i) + 2 k d_{z_\al}
  G_{2\la} (z_\al,P_i) ] \res_{z = P_i} f(\wt m)_\la
  (z,(z_\beta)_{\beta\neq \al}) ,
\end{align}
where the functions $G$ and $G_{2\la}$ are (twisted) Green functions.

The relation to the usual formulation of the KZ connection in the
rational case is the following. In that case, the KZ connection has
the form
\begin{equation} \label{KZ:class}
  2(k+2)\nabla_{P_i}\psi(P_i) = 2(k+2)\pa_{P_i}\psi((P_i)_i) - K^{rat}_i
  \psi(P_i),
\end{equation} 
with $\psi(P_i)$ in a tensor product $\otimes_i V_{\La_i}$ of lowest
weight $\bar\G$-modules, and $K^{rat}_i = \sum_{j\neq i}
{{t^{(ij)}}\over{P_i - P_j}}$. Equation (\ref{cite:KZ}) above may be
viewed as the expression of the action of $K_i$ on ``Bethe ansatz
vectors'' $\wt e(\zeta_1) \cdots \wt e(\zeta_k)(\otimes_i
v^{bot}_{\La_i})$, where $\wt e(z) = \sum_i {{e^{(i)}}\over{z -
    P_i}}$. Extracting coefficients of $\prod_{} (\zeta_i -
z_j)^{a_{ij}}$ from (\ref{cite:KZ}), one recovers (\ref{KZ:class}).
The equation for the bottom component of $\psi(P_i)$ is simpler than
(\ref{cite:KZ}) (see eq.  (\ref{var:pts:0})).

The operators $(T_z)_{z\in X}$ depend in a simple way on the level
$k$. In sect.  \ref{comm:do}, we show that these operators commute
when $k$ is critical, thus defining a commuting family of differential
operators, acting on a finite-dimensional bundle over the degree zero
part $J^0(X)$ of the Jacobian of $X$ (Thm.\ \ref{main}). This is proved
using a class of modules $W_{n|m,m'}$ generalizing the twisted Weyl
modules. 

In the case where there are no $z^{(i)}_j$, these operators take the
form
$$ T_z f(\la) = \left( \sum_\nu (\sum_a \omega_a(z) \pa_{(h_\nu)_a})^2
  + \sum_{\al \in \Delta_+} \sum_a D_z^{(\la,\al)}
  \omega_a(z)\pa_{(\al^\vee)_a} + k\sum_{\al\in \Delta_+}
  \omega_{(\la,\al)}(z) \right) f(\la),
$$ $(h_\nu)$ an orthonormal basis of the Cartan subalgebra $\bar\HH$
of $\bar\G$, $\omega_a$ the canonical one-forms on $X$, $\Delta_+$ the
set of positive roots of $\bar \G$, $\la$ is a collection
$(\la_1,\ldots,\la_g)$ of variables in $\bar\HH$, $\al^\vee$ is the
coroot associated to the root $\al$, $D_z^{(\la^{(i)})}$ is a
connection depending on $(\la^{(i)})$ in $\CC^g$, on the canonical
bundle $\Omega_X$ (see (\ref{D:la})) and $\omega_{(\la^{(i)})}$ is a
quadratic differential form depending on the same variables (see
(\ref{omega:la})).

We close the paper by explaining the link of the operators
$(T_z)_{z\in X}$ with the Beilinson-Drinfeld (BD) operators (Rem.
\ref{Beil-Dr}).

In a sequel to this paper, we will construct $q$-deformations of the
operators $T_z$, by replacing the inclusion $U\G^{out} \subset U\G$ by
some inclusion of quasi-Hopf algebras, which were introduced in work
of one of us and V. Rubtsov.  The outcome will be a commuting family
of difference-evaluation operators, which may be viewed in the
rational case as the Bethe ansatz formulation of the qKZ operators.

One may hope to obtain hypergeometric representation for solutions of
the KZB equations formulated in sect. \ref{KZB}. This may be related
with the formulas of \cite{Gaw} expressing the scalar product on the
space of conformal blocks. 

We would like to thank B. Feigin, V. Rubtsov and V. Tarasov for
discussions on this paper. B.E. would like to thank the FIM, ETHZ, for
hospitality at the time this work was done.

\section{Bases of functions on $X$} \label{bases}

Let $X$ be a smooth, compact complex curve; denote by $g$ its genus.
Let $P_0$ be a point of $X$. Denote by $\cK$ and $\cO$ the local field and
ring of $X$ at $P_0$.  Denote by $\Omega_{\cK}$ and $\Omega_{\cO}$ the
spaces of differentials and regular differentials at the formal
neighborhood of $P_0$. The residue defines a natural pairing between
$\cK$ and $\Omega_{\cK}$.

In what follows, we will fix a system $(A_a,B_a)_{a=1,\ldots,g}$ of
$a$- and $b$-cycles on $X$. We will denote by $\gamma_{A_a}$ and
$\gamma_{B_a}$ the corresponding deck transformations on the universal
cover $\wt X$ of $X$, and by $\sigma$ the projection from $\wt X$ to
$X$.

Define $R_{(b)}$ as the set of functions defined on $\wt X$, regular
outside $\si^{-1}(P_0)$, such that $f(\gamma_{A_a}z) = f(z)$ and
$f(\gamma_{B_a}z) =f(z) + \al_a(f)$, with $\al_a(f)$ some constants.
Let us also denote by $R$ the space of functions on $X$, regular
outside $P_0$.

\begin{prop} \label{codim}
  $R_{(b)} \cap\cO = \CC 1$. $R$ has codimension $g$ in $R_{(b)}$. Moreover,
  $R_{(b)} + \cO = \cK$.
\end{prop}

{\em Proof.}
The first point is clear: for any $f$ in $R_{(b)} \cap \cO$, $df$ is
a regular form with vanishing $a$-periods, and therefore vanishes. 
  
To prove the second point, define $R_{(ab)}$ as the set of regular
functions defined on the universal cover of $X-P_0$, such that
$f(\gamma_{A_a}z) = f(z) + \beta_a(f)$ and $f(\gamma_{B_a}z) = f(z) +
\al_a(f)$, with $\al_a(f)$ and $\beta_a(f)$ some constants.  We will
show that $R$ has codimension $2g$ in $R_{(ab)}$.  $R_{(ab)} \cap \cO$
has dimension $g+1$ (it is spanned by the constants and the
$\int_{P_0}^x \omega_a$). On the other hand, we have $R_{(ab)} + \cO =
\cK$, because $\cK/R_{(ab)} + \cO$ is zero (the differential maps it
injectively to $\Ker \res / \Omega_R + \Omega_{\cO}$ where $\res$ is
the residue map from $\Omega_\cK$ to $\CC$, which is the kernel of the
residue map from $H^1(X,\Omega_X)$ to $\CC$ and is therefore zero). We
have an exact sequence $0 \to (R_{(ab)}\cap \cO)/ (R\cap \cO) \to
R_{(ab)}/R \to (R_{(ab)}+\cO)/(R+\cO) \to 0$, therefore $\dimm(R_{(ab)}/R) =
\dimm(R_{(ab)}\cap \cO / R\cap \cO) + \dimm(\cK / R + \cO) = 
2g$. Since $\dimm(R_{(ab)}/R) = \dimm(R_{(ab)}/R_b)+ \dimm(R_{(b)}/R)$, we
have $\dimm(R_{(ab)}/R_{(b)})+ \dimm(R_{(b)}/R) = 2g$.

On the other hand, $\dimm(R_{(ab)}/R_{(b)})$ and $\dimm(R_{(b)}/R)$
are both $\leq g$, because the maps $R_{(ab)}/R_{(b)} \to \CC^g$
sending the class of $f$ to $(\beta_a(f))$ and $R_{(b)}/R \to \CC^g$
sending $f$ to $(\al_a(f))$, are both injections.

It follows that $\dimm(R_{(ab)}/R_{(b)})$ and $\dimm(R_{(b)}/R)$ are
both equal to $g$.

Finally, the fact that $\cK / \cO + R$ is equal to $H^1(X,\cO_X)$ and
has therefore dimension $g$ implies the last point.  
\hfill \qed
\medskip

\begin{cor} \label{r_a}
  For $a=1,\ldots,g$, we have functions $r_a$ defined on $\wt X$,
  regular outside $\si^{-1}(P_0)$, with the properties 
  $$ r_a(\gamma_{A_b}z) = r_a(z), \quad r_a(\gamma_{B_b}z) = r_a(z) -
  \delta_{ab},
  $$ for $b = 1, \ldots,g$. The functions $r_a$ are well-defined up to
  addition of functions of $R$.
\end{cor}

Fix a coordinate $z$ at $P_0$.  Let us denote by $\mm$ the maximal
ideal of $\cO$, by $(r_{i,0}^{in})$ a basis of $\mm$ and by
$(r_i^{out},1)$ a basis of $R = H^0(X-P_0,\cO_X)$, such that
$\res_{P_0}r^{out}_i {{dz}\over z} = 0$.  From Prop.\ \ref{codim}
follows that we can fix functions $(r_{a})_{a = 1, \ldots,g}$ of $R_{(b)}$
such that $\res_{P_0}r_a {{dz}\over z} = 0$, so that
$(r_a,r_i^{out},1)$ is a basis of $R_{(b)}$ and
$(r_{i,0}^{in},r_a,r^{out}_i,1)$ is a basis of $\cK$.

Let $(\omega_a)_{a=1,\ldots,g}$ be the basis of the space of
holomorphic differentials $\Omega_{\cO} \cap H^0(X-P_0,\Omega_X)$, dual to
$(r_{a})$.  We have 
$$
{1\over{2i\pi}}\int_{A_a}\omega_b = \delta_{ab}. 
$$
We can fix families $(\omega^{in}_i)$ and
$(\omega^{out}_i)$ in $\Omega_{\cO}$ and $H^0(X-P_0,\Omega_X)$, so that
$(\omega_i^{out},\omega_a,\omega^{in}_i,{{dz}\over z})$ is the basis of
$\Omega_{\cK}$ dual to $(r_{i,0}^{in},r_a,r^{out}_i,1)$.

We associate with these dual bases the Green function defined as 
\begin{equation} \label{green:0}
G(z,w) = \sum_{i}\omega_{i}^{out}(z) r_{i,0}^{in}(w). 
\end{equation}
It is clear that $G$ depends only on the choice of $a$-cycles in $X$.

Denote by $J(X)$ the Jacobian of $X$. It is the direct sum of its
degree $n$ components $J^{n}(X)$, with $n$ integer, which are
identified with the sets of classes of line bundles of degree $n$ on
$X$.  Denote by $\Gamma$ the lattice of periods of $X$, which we
identify with a lattice in $\CC^g$ via the basis dual to
$(\omega_a)_{a = 1 ,\ldots,g}$.  $J^{0}(X)$ is identified with the
quotient $\CC^g/\Gamma$, as follows: for some $\la = (\la_a)$ in
$\CC^g$, the corresponding line bundle is denoted by $\cL_\la$.
Sections of $\cL_{\la}$, regular outside a finite subset $S$ of $X$,
are identified with the functions on the universal cover of $X$,
regular outside the preimage of $S$, such that $f(\gamma_{A_a}z) =
f(z)$ and $f(\gamma_{B_a}z) = e^{\la_a}f(z)$.  Multiplication by the
functions $e^{\int^z \omega_a}$ identifies the spaces of sections of
$\cL_\la$ and $\cL_{\la'}$, for $\la$ and $\la'$ in the same class of
$\CC^g/\Gamma$.

In what follows, we will set 
\begin{equation} \label{not}
R_\la = H^0(X-\{P_0\} , \cL_\la) . 
\end{equation} 
Let $\la$ be nonzero in $J^{0}(X)$. We may identify $H^0(X-\{P_0\} ,
\Omega_X\otimes \cL_\la)$ with the space of $1$-forms $\omega$ on
the universal cover of $X$, regular outside the preimage of $P_0$,
such that $\gamma_{A_a}^*(\omega) = \omega$ and
$\gamma_{B_a}^*(\omega) = e^{\la_i}\omega$.  The space
$H^0(X,\Omega_X \otimes \cL_\la)$ may be identified with the
intersection $\Omega_{\cO} \cap H^0(X-\{P_0\} , \Omega_X \otimes
\cL_\la)$.  By the Riemann-Roch theorem, it has dimension $g-1$.
Let $(\omega_{a;\la})_{a = 1,\ldots,g-1}$ be a basis of this space.
We may complete it to a basis $(\omega_{i;\la}^{out} ,
\omega_{a;\la} , \omega_{i}^{in})$ of $\Omega_\cK$, such that
$(\omega_{i;\la}^{out} , \omega_{a;\la})$ is a basis of
$H^0(X-\{P_0\} , \Omega_X \otimes \cL_\la)$ and $(\omega_{a;\la}
, \omega_{i}^{in})$ is a basis of $\Omega_\cO$.  Moreover, we may
assume that the $\omega_i^{in}$ have a zero of order $\ge g-1$ ar
$P_0$ (for example, we may choose $\omega_i^{in} = z^{g-1+i}dz$, $i\ge
0$).

Let $(r_i^{in},r_{a;-\la},r^{out}_{i;-\la})$ be the basis of $\cK$ dual
to $(\omega_{i;\la}^{out} , \omega_{a;\la} , \omega_{i}^{in})$.  Then
$(r^{in}_{i})$ is a basis of $\cO$ and $(r^{out}_{i;-\la})$ is a basis
of $H^0(X-\{P_0\} , \cL_\la^{-1})$. The assumption on zeroes of the
$\omega_i^{in}$ implies that the $r_{a;-\la}$ have poles at $P_0$ of
order $\leq g-1$.

The twisted Green function defined by these bases is
\begin{equation} \label{green:la}
G_\la(z,w) = \sum_{a=1}^{g-1} \omega_{a;\la}(z)r_{a;-\la}(w)
+  \sum_{i} \omega_{i;\la}^{out}(z) r^{in}_i(w). 
\end{equation}

\begin{remark} 
  {\it Expression of the Green functions.} We may set
  $$ r_a(z) = \pa_{\epsilon_a}\ln \Theta(-A(z) + g A(P_0)-\Delta) ,
  $$ where $\Theta$ is the Riemann theta-function on $J^{0}(X)$,
  $\Delta \in J^{g-1}(X)$ is the vector of Riemann constants of $X$,
  $\epsilon_a$ is the $a$-th basis vector of $\CC^g$ and $A$ is the Abel
  map from $X$ to $J^{1}(X)$.

  A formula for $G_{\la}$ is
\begin{align*} 
  G_\la(z,w) & = {{\Theta(A(z)-A(w) + (g-1) A(P_0) - \la -
      \Delta)}\over{\Theta(A(z)-A(w)+ (g-1) A(P_0) -
      \Delta)\Theta((g-1)A(P_0) - \la - \Delta)}} \cdot \\ & \cdot
  \sum_{i=1}^g {{\pa \Theta}\over{\pa \la_a}}((g-1)A(P_0) - \Delta)
  \omega_a(z) ;
  \end{align*} 
  $G_{\la}(z,w)$ is a $\la$-twisted $1$-form in $z$, with simple pole
  at $z = w$ and residue $1$, and $g-1$ zeroes at $z = P_0$; it is
  also a
  $(-\la)$-twisted function in $w$, with simple poles at $w = z$ and a
  pole of order $g-1$ at $w = P_0$. This is because $\sum_{i=1}^g
  {{\pa \Theta}\over{\pa \la_a}}((g-1)P_0-\Delta) \omega_a(z)$, which
  is equal to $-d_z \Theta(w-z+(g-1)P_0-\Delta)_{|w=z}$, is a
  holomorphic form with $g-1$ zeroes at $P_0$. For $z,w$ fixed,
  $G_\la(z,w)$ is a meromorphic function in $P_0$. One may replace
  $(g-1)P_0$ by any effective divisor $Q = \sum_i n_i Q_i$ of degree
  $g-1$ in the definition of $G_\la$, and obtain this way
  $G_\la^Q(z,w)$, a $\la$-twisted $1$-form in $z$, with simple pole at
  $w$ and $n_i$ zeroes at the $Q_i$, which is also a $(-\la)$-twisted
  function in $w$, with simple poles at $z$ and poles of order $n_i$
  at $w = Q_i$, and is a meromorphic function in the $Q_i$.

  A formula for $G(z,w)$ is 
  $$ G(z,w) = d_z \ln \Theta(A(w)-A(z)+(g-1)A(P_0)-\Delta) - d_z
  \ln\Theta(gA(P_0)-A(z) -\Delta).
  $$ $G(z,w)$ is a one-form in $z$ with simple pole at $w$ and residue
  $1$; simple pole at $z=P_0$, regular at other points, and such that
  $\int_{A_a} G( \cdot ,w) = 0$ for $w$ near $P_0$; and a function in $w$,
  multivalued in $w$ around $b$-cycles, such that $G(z,\gamma_{B_a}w)
  = G(z,w) + \omega_a(z)$, vanishing for $w = P_0$, with simple pole
  at $w=z$, and regular at other points.

  These properties of $P_0$ imply that two $G(z,w)$ attached to
  different points $P_0$ differ by a form in $z$, constant in $w$. In
  what follows, we will set
  \begin{equation} \label{wt:omega}
    \wt \omega(z,w) = d_w G(z,w).
  \end{equation} 
  $\wt \omega(z,w)$ is a bidifferential form in $z,w$ with the local
  expansion at any point of $X$, $\wt \omega(z,w) =
  {{dzdw}\over{(z-w)^2}} + r(z) dzdw + O(z-w)dzdw$.

  $\wt\omega$ is symmetric in $z$ and $w$, because $\wt\omega(z,w) -
  \wt\omega(w,z)$ has no poles and for $w$ near $P_0$, $\int_{A_a}
  \wt\omega(\cdot,w) - \wt\omega(w,\cdot) = d_w \int_{A_a} G(\cdot,w)
  - (G(w,\gamma_{A_a}z) - G(w,z)) = 0$ because $\int_{A_a}G(\cdot,w) =
  0$ and because $G(w,\cdot)$ is single-valued along $a$-cycles.  That
  $\wt \omega$ is symmetric can also be viewed as a consequence of the
  expression $\wt\omega = d_z d_w \ln \Theta(A(w)-A(z) +
  \delta-\Delta)$ where $\delta$ in $J^{g-1}(X)$ is some odd
  theta-divisor.
\end{remark}

\section{Twisted conformal blocks}

\subsection{Twisted conformal blocks} \label{int}

Let $\bar\G$ be a simple complex Lie algebra. Let us set $\G = (\bar\G
\otimes \cK) \oplus \CC K$, $\G^{in} = (\bar\G \otimes \cO) \oplus \CC
K$, $\G^{out} = \bar\G \otimes R$. For $x$ in $\bar \G$,
$\eps$ in $\cK$, we set $x[\eps] = (x \otimes \eps,0)$; the commutation
rules on $\G$ are then
$$
[x[\eps],y[\eps']] = [x,y][\eps\eps'] + K\langle d \eps, \eps'\rangle 
(x|y),  
$$ with $(\cdot |\cdot )$ the invariant scalar product on $\bar\G$
such that $(\theta^\vee|\theta^\vee)=2$, where $\theta^\vee$ is the
coroot associated to a maximal root $\theta$, and $\langle
\omega,\eps \rangle = \res_{P_0}(\omega\eps)$. We view $\G^{out}$ as a
subalgebra of $\G$, using the embedding $x\otimes r\mapsto x[r]$.

Let $V$ be a $\G$-module of level $k$, and let $\psi$ be a
$\G^{out}$-invariant linear form on $V$.  Fix a Cartan decomposition
$\bar\G = \bar \HH \oplus \bar\N_+ \oplus \bar \N_-$. Let $r$ be the
rank of $\bar\G$. Let $\Delta$ be the set of roots of $\bar\G$, and
define the positive roots as those associated with $\bar\N_+$. For
each $\al$ in $\Delta$, define $\bar\G_\al$ as the root subspace of
$\bar\G$ associated with $\al$.  For each simple root $\al_i$, let us
fix $e_i,h_i$ and $f_i$ in $\bar\G_{\al_i}$, $\bar\HH$ and
$\bar\G_{-\al_i}$, such that $(e_i,h_i,f_i)$ is an $\SL_2$-triple.

Let $(r_a)_{a = 1 , \ldots,g}$ be as in Cor.\ \ref{r_a} and let
$(\la^{(i)}_a)_{a = 1,\ldots,g, i = 1, \ldots,r}$ be formal variables
and define the linear form $\psi_\la$ in $V$
\begin{equation} \label{formal:psi:la}
  \langle \psi_\la , v \rangle = \langle \psi ,
  e^{\sum_{i,a} \la^{(i)}_a h_i[r_a]} v \rangle.
\end{equation} 
This form is independent of the choice of the $r_a$, because
$[h_i[r],h_j[r_a]] = 0$ for $r$ in $R$. 

In the case where $V$ is an integrable module, one expects that one
can make sense of (\ref{formal:psi:la}) for complex $\la$. If one
wished to argue that the action of $\G$ on $V$ lifts to a projective
action of the associated Kac-Moody group, one would meet the
difficulty that the functions $e^{\sum_a \la_a^{(i)} r_a}$ have
essential singularities at $P_0$, so that we cannot view
$e^{\sum_{a,i} \la^{(i)}_a h_i[r_a]}$ as an element of the Kac-Moody
group.

However, we have:

\begin{thm} \label{seerose} 
  For $\psi$ a $\G^{out}$-invariant form on $L_{\La,k}$, the form
  $\psi_\la = \psi \circ e^{\sum_{i,a} \la^{(i)}_a h_i[r_a]} $
  on $L_{\La,k}$ has the following properties:

  1) For any $v$ in $L_{\La,k}$, the function $\langle \psi_\la ,
  v\rangle$ is the formal expansion at $0$ of an analytic function in
  $\la$, which satisfies the equations
  $$ \pa_{\la^{(i)}_a} \langle \psi_\la , v\rangle = \langle \psi_\la ,
  h_i[r_a]v\rangle,
  $$ $a = 1, \ldots,g, i = 1, \ldots, r$. 

  2) Set $\la_a = \sum_i \la^{(i)}_a h_i$. Set $\la =
  (\la_1,\ldots,\la_g)$ and $$ \G_\la^{out} = (\bar\HH \otimes R)
  \oplus \oplus_{\al\in \Delta} (\bar\G_\al \otimes R_{\langle
    \al,\la_1 \rangle , \ldots , \langle \al,\la_g \rangle }).
  $$ Then $\psi_\la$ is a $\G^{out}_\la$-invariant form on
  $L_{\La,k}$.

  3) For any $v$ in $L_{\La,k}$, the function $\la \mapsto \langle
  \psi_\la , v\rangle$ has the following theta-like behavior. Set
  $\omega_{ab} = \int_{B_b}\omega_a$, $\zeta_a(z) = \int_{P_0}^z
  \omega_a$, and $\Omega_a = \sum_b \omega_{ab}\delta_b$, with
  $\delta_a$ the $a$-th basis vector of $\CC^g$. Then
  $$ \langle \psi_{\la^{(1)} , \ldots, \la^{(l)}+2i\pi\delta_a,
    \ldots, \la^{(r)}} , v \rangle = \langle
  \psi_{\la^{(1)},\ldots,\la^{(r)}} , v \rangle$$ and
  $$ \langle \psi_{\la^{(1)} , \ldots, \la^{(l)}+ 2i\pi\Omega_a,
    \ldots, \la^{(r)}} , v \rangle = e^{- k (h_l|\la_a)- {1\over
      2}i\pi k \omega_{aa} (h_l|h_l)}\langle \psi_{\la^{(1)} , \ldots,
    \la^{(r)}} , e^{h_l[\zeta_a]}v \rangle,
  $$ where $\la_a = \sum_{i=1}^r \la_a^{(i)}h_i$.
\end{thm}

{\em Proof.} See the appendix. 

\subsection{Twisted correlation functions in the $\SL_2$ case}

In this section, we assume that $\bar\G = \SL_2$ and $\La = 0$.  Let
$\psi$ be a $\G^{out}$-invariant form on $L_{0,k}$. Let $z$ be a local
coordinate at $P_0$ and set $e(w) = \sum_{i\in\ZZ}e[z^i]w^{-i-1}dw$.
For $n$ a positive integer, set $n = a k + b$, $0\leq b < k$, and
$v_{n} = f[z^{-2a-1}]^{b}v_{[a]}$, with $v_{[a]}$ as in Lemma
\ref{5.6}. We have $h[1]v_n = -2n v_n$.  Set $f_\la(z_1,\ldots,z_n) =
\langle \psi_\la, e(z_1) dz_1\cdots e(z_n) dz_n v_n \rangle$.

\begin{prop} \label{funct} (see \cite{FS})
  The form $f_\la(z_1,\ldots,z_n)$ depends analytically on $\la$
  in $J^0(X)$ and the $z_i$ in $X - \{P_0\}$. It satisfies the relations 
  \begin{equation} \label{int1} f_{\la + 2i\pi \delta_a}(z_1,\ldots,z_n) =
    f_{\la}(z_1,\ldots,z_n),
  \end{equation} and
  \begin{equation} \label{int2} 
    f_{\la + 2i\pi \Omega_a}(z_1,\ldots,z_n) = e^{-k(h|h)\la_a -
      {1\over 2}i\pi k \omega_{aa} (h|h)} e^{2\sum_{l=1}^n
      \int_{P_0}^{z_l} \omega_a} f_{\la}(z_1,\ldots,z_n) .
  \end{equation} 
  Moreover, it depends on $z_i$ as a section of $\Omega_X\cL_{-2\la}$,
  regular on $X$ except for a pole of order $\leq 2a+2-2\delta_{b,0} $
  at $P_0$; it is symmetric in the $z_i$, and vanishes if $k+1$
  variables $z_i$ coincide.
\end{prop}

{\em Proof.} The proof is analogous to that of \cite{FS}. Identities
(\ref{int1}) and (\ref{int2}) follow from Thm.\ \ref{seerose}, 3) and
the commutation relation $[h[\zeta_a],e(z_i)] = 2 \zeta_a(z_i)e(z_i)$.
\hfill \qed \medskip

Since $(h|h) = 2$, we have $$f_{\la}(z_1,\ldots,z_n) = \sum_{l =
  1}^{(2k)^g} \Theta_{2k}^{[l]}(\la + {1\over k}\sum_{i=1}^n A(z_i))
f^{[l]}(z_1,\cdots,z_n),$$ where the $\Theta_{2k}^{[l]}$ are a basis
of the space of $2k$th order theta-functions on $J^0(X)$.

\begin{remark} \label{vopros}
  With $f_{FS}(z_1,\ldots,z_n)$ the forms introduced in \cite{FS},
  $f_{0}(z_1,\ldots,z_n)$ coincides with $f_{FS}(z_1,\ldots,z_n)$.  It
  is not clear what are the functional properties of the
  $f^{[l]}(z_1,\cdots,z_n)$, and how to obtain the
  $f^{[l]}(z_1,\cdots,z_n)$ directly from $f_{FS}(z_1,\ldots,z_n)$.
\end{remark}

\begin{remark}
  The forms $f_{\la}(z_1,\ldots,z_n)$ provided by conformal blocks
  also satisfy some vanishing conditions at $\la=0$ (see \cite{FW}).
  These conditions, together with the functional properties of Prop.\
  \ref{funct}, should probably characterize these forms.
\end{remark}

\section{Lifts of generalized theta-functions to $Bun_B$} \label{lifts}

{}From the works \cite{BL,Kumar} follows that conformal blocks may be
viewed as the space of sections of a line bundle on the moduli space
$Bun_{\bar G}$ of principal $\bar G$-bundles over an complex curve
$X$, for $\bar G$ the simply connected group associated with $\bar\G$.
This identification is as follows: $Bun_{\bar G}$ is identified with
the double coset $\bar G(R)\backslash \bar G(\cK) / \bar G(\cO)$, with
$\cK$ the local field at some point $P_0$ of $X$, $\cO$ the local ring
at $P_0$ and $R$ the ring of functions regular outside $P_0$. For $k$
integer $\ge 0$, the level $k$ vacuum representation $L_{0,k}$ of the
Kac-Moody algebra $\G$ associated with $\bar \G$ carries a projective
representation of $\bar G(\cK)$. Fix a lift $x \mapsto \wt x$ of $\bar
G(\cK)$ to its universal central extension. Let $\G^{out}$ be the Lie
algebra $\bar\G \otimes R$. To each $\G^{out}$-invariant form
$\psi^{out}$ on $L_{0,k}$ is associated the function
\begin{equation} \label{al}
g \mapsto \langle \psi^{out},\wt g v_{top}\rangle
\end{equation}
on $\bar G(\cK)$, where $v_{top}$ is the vacuum vector of $L_{0,k}$,
which is a section of a power of the determinant bundle over
$Bun_{\bar G}$. This construction can be extended to the case of
marked points and integrable representations other than $L_{0,k}$.  In
what follows, we will consider the situation of some integrable module
$L_{\La,k}$ at $P_0$, with highest weight vector $v_{top}^{(P_0)}$.

It was proposed to study these functions through their lifts to moduli
spaces of flags of bundles (\cite{Bertram,Thadd}).  In \cite{FS},
Feigin and Stoyanovsky studied the lift of conformal blocks to a
space, which in the case $\bar\G = \SL_n$ can be described as
$Bun_{(n_i,P_0)}$, the moduli space of bundles with filtration
$E_1\subset E_2 \subset \cdots $ and associated graded isomorphic to
$\oplus_i \cO(n_i P_0)$, $n_i$ some integer numbers.  Since this space
is isomorphic to $N(R)\backslash N(\cK)diag(z^{n_i})/N(\cO)$, with $N$
the maximal unipotent subgroup of $\bar G$, lifts of functions
provided by the conformal blocks are the
\begin{equation} \label{n:corr}
\langle \psi^{out} , n_{\cK} (w v_{top}^{(P_0)}) \rangle,
\end{equation} 
$n_{\cK}$ in $N(\cK)$, $w = diag(z^{n_i})$ an affine Weyl group
translation. Generating functions for these quantities are the forms
$$\langle \psi^{out} , \prod_{i\ \on{ simple}}
\prod_{j=1}^{n_j}e_i(z^{(i)}_j) dz^{(i)}_j (w v_{top}^{(P_0)})\rangle,
$$ where $e_i(z)dz$ are the currents associated to the nilpotent
generators $e_i$ attached to the simple roots of $\bar\G$.  In
\cite{FS}, Feigin and Stoyanovsky characterized the functional
properties of these forms.

Let us study the lift of functions (\ref{al}) to $Bun_B$, the moduli
space of $B$-bundles over $X$, where $B$ is the Borel subgroup of
$\bar G$ containing $N$. $Bun_B$ can be described as the double
quotient $B(K)\backslash B(\AAA)/B(\cO_\AAA)$, where $K$ is the
function field $\CC(X)$, $\AAA$ is the adeles ring of $X$ and
$\cO_\AAA$ its subring of integral adeles. To make sense of the
analogue of (\ref{n:corr}) for the space of $B$-bundles, one should
replace the representation at $P_0$ by its ``adelic'' version
$L^\AAA$, which is its restricted tensor product with vacuum
representations at the points of $X - \{P_0\}$. To $\psi^{out}$ is
then associated a $\bar\G \otimes K$-invariant form $\psi^\AAA$ (see
Lemma \ref{finkel}).  In the case of $B$-bundles, lifts of the
functions on $Bun_{\bar G}$ provided by conformal blocks are the
\begin{equation} \label{b:corr} 
  b \mapsto \langle \psi^{\AAA}, b v^\AAA_{top} \rangle ,
\end{equation} 
for $b \in B(\AAA)$, $v^\AAA_{top}$ the product of the highest weight
vector of the module at $P_0$ with the vacuum vectors at other points.
$b$ can be decomposed as a product $ntw$, with $n$ in $N(\AAA)$, $t$
in $T(\AAA)$ with all components of degree zero ($T$ is the Cartan
subgroup associated to $B$; the degree in $\AAA^\times$ is defined as
the sum of the valuations of all components) and $w$ a product of
affine Weyl group translations.  In the case $\bar\G = \SL_n$, $b$
represents a filtered bundle whose associated graded is a sum of line
bundles, associated to the projections in the Jacobian $J(X) =
K^\times \backslash \AAA^\times / \cO_\AAA^\times$ of the components
of $tw$.

The computation of (\ref{b:corr}) may be done as follows. $w
v_{top}^\AAA$ is an extremal vector of $L^\AAA$. $n$ may be replaced
by an element $n_{\cK}$ of $N(\AAA)$ with only nontrivial component at
$P_0$. The map $\la \mapsto f_\la$ of sect. \ref{tame} is a section of
the projection map $K^\times \backslash (\AAA^\times)^0 \to J^0(X)$
(the $^0$ denotes the degree zero parts). $t$ can be decomposed as
$t^{out}t_\la t^{in}$, $t^{out}$ in $T(K)$, $t^{in}$ in $T(\cO_\AAA)$
and $t_\la = \prod_i t_i[f_{\la^{(i)}}]$, $t_i$ the subgroups of $\bar
G$ associated to the simple coroots of $\bar\G$. Then (\ref{b:corr})
is equal to $(t^{in},n_\cK) \langle \psi^\AAA , t_\la n_\cK (w
v_{top}^\AAA)\rangle$ (where $(,)$ denotes the group commutator).

Therefore to compute (\ref{b:corr}), it suffices to compute the
\begin{equation} 
  \langle \psi^{\AAA}, \prod_{i=1}^r t_i[f_{\la^{(i)}}] \prod_{i=1}^r
  e_i[\eps^{(i)}_1] \cdots e_i[\eps^{(i)}_{n_i}] (wv_{top}^\AAA)
  \rangle, 
\end{equation}
where $r$ is the rank of $\bar\G$.  In Thm.\ \ref{seerose}, we study the
linear form
\begin{equation} \label{wl}
  v\mapsto \langle \psi^\AAA, \prod_{i=1}^r t_i[f_{\la^{(i)}}] (v
  \otimes \otimes_{x\in X-\{P_0\}} v_{top}) \rangle,
\end{equation}
for $v$ in $L_{\La,k}$. 

{}From Thm.\ \ref{seerose} follows that the expansion at $(\la_a^{(i)})
= 0$ of (\ref{tw:corr}) is equal (up to multiplication by phase
factor) to
\begin{equation} \label{tw:corr}
  \langle \psi^{out}, e^{\sum_{i,a} \la_a^{(i)} h_i[r_a]}
  \prod_{i}e_i[\eps^{(i)}_1] \cdots e_i[\eps^{(i)}_{n_i}] (w
  v_{top}^{(P_0)}) \rangle. 
\end{equation} 
Generating functions for (\ref{tw:corr}) are the forms (\ref{oo}).

The interest of expressing (\ref{b:corr}) in the form (\ref{tw:corr})
is that the latter expression is computed in a single module located
at $P_0$.  When the $\la_a^{(i)}$ are formal, (\ref{tw:corr}) also
makes sense in arbitrary modules. What we will do now is compute the
action of the Sugawara tensor on these correlation functions. 


\section{Expression of the KZB connection}

\subsection{Action of the Sugawara tensor on the twisted correlation functions 
  ($\bar\G = \SL_2$)}

In this section, we treat the case $\bar\G = \SL_2$. Let $n$ be an
integer and let $v_n$ be a vector of $L_{\La,k}$ such that $h[1]v_n =
-2n v_n$, $h[t^{k}]v_n = 0$ for $k > 0$ and $f[t^{k}]v_n = 0$ for
$k\geq -(g-1)$. An example of $v_n$ is in the vacuum module $L_{0,k}$,
the extremal vector $f[t^{-(2a-1)}]^k \cdots f[z^{-1}]^k v_{top}$,
with $2a+1 \geq g-1$.

In what follows, we will denote by $x(z)$ the series $\sum_{i\in\ZZ}
x[t^i]z^{-i-1}dz$, for $x$ in $\bar\G$.

The expression for the Sugawara tensor is
\begin{align} \label{expr:T}
  & 2(k+2)T_{\wt \omega}(z) \\ & \nonumber = \lim_{z\to z'}[e(z) f(z')
  + f(z) e(z') + {1\over 2} h(z) h(z') -3k\wt \omega(z,z')],
  \end{align}
  with $\wt\omega$ as in (\ref{wt:omega}). It is used to define the
  KZB connection in sect. \ref{KZB}.

\subsubsection{Action of the currents on the correlation functions}

Assume that $m$ is $\le -(g-1)$. Let us compute some correlation
functions in $L_{\La,k}$. 

\begin{lemma} \label{lemma:act:h}
  We have
  \begin{align*} & \langle \psi_{\la}, h(z) e(z_1)\cdots e(z_n) v_n \rangle
    \\ & = \left( \sum_a \omega_a(z) \pa_{\la_a} + 2
      \sum_{\al=1}^n G(z,z_{\al})\right) f_{\la}(z_1 , \cdots , z_n),
  \end{align*} where $G(z,z_\al)$ is as in (\ref{green:0}). 
\end{lemma}

{\em Proof.} Let us write $h(z)  = \sum_i h[r^{out}_i] \omega_i^{in}
+ \sum_a h[r_a]\omega_a + \sum_i h[r^{in}_{i,0}]\omega_i^{out}$. The
contribution of the first term of this sum is zero by invariance of
$\psi_\la$, the contribution of the second part is the differential
part.  The contribution of the third part is
\begin{align*}
  & \sum_i \langle \psi_\la, h[r^{in}_{i,0}] e(z_1) \cdots e(z_n) v_n
  \rangle \omega_i^{out}(z) \\ & = \sum_i \sum_{j=1}^n 2
  r_{i,0}^{in}(z_j)\langle \psi_\la, e(z_1) \cdots e(z_n) v_n \rangle
  \omega_i^{out}(z) + \sum_i \langle \psi_\la, e(z_1) \cdots e(z_n)
  h[r^{in}_{i,0}] v_n \rangle \omega_i^{out}(z) \\ & = \sum_i \sum_{j=1}^n
  2 G(z,z_j) \langle \psi_\la, e(z_1) \cdots e(z_n) v_n \rangle,
\end{align*}
because $v_n$ is annihilated by the positive Cartan modes. \hfill \qed
\medskip

\begin{lemma} \label{lemma:act:f}
  We have 
\begin{align} \label{act:f}
  & \langle \psi_{\la}, f(z) e(z_1)\cdots e(z_{n+1}) v_{n} \rangle = \\ &
  \nonumber - \sum_{\al} G_{2\la}(z,z_\al) \left(
    \sum_{a} \omega_{a}(z_\al)\pa_{\la_a} + 2 \sum_{\beta\neq \al}
    G(z_\al, z_{\beta}) \right) \\ & \nonumber \langle
  \psi_{\la}, e(z_1)\cdots e(z_{\al-1})e(z_{\al+1}) \cdots e(z_{n+1})
  v_{n} \rangle \\ & \nonumber + k \sum_{\al=1}^{n+1}
  d_{z_{\al}}G_{2\la}(z,z_{\al}) \langle \psi_{\la} , e(z_1)\cdots
  e(z_{\al-1}) e(z_{\al+1}) \cdots e(z_{n+1}) v_{n}\rangle, 
\end{align}
with $G_{2\la}(z,z_\al)$ as in (\ref{green:la}).
\end{lemma}

{\em Proof.} Write
$$ f(z) = \sum_i f[r_{i;-2\la}^{out}] \omega_i^{in}(z) + \sum_a
f[r_{a;-2\la}] \omega_{a;2\la}(z) + \sum_i f[r_i^{in}]
\omega_{i;2\la}^{out}(z). 
$$ 
The contribution of the first term is zero by invariance of $\psi_\la$. 
The contributions of the two next terms is 
\begin{align} \label{interm}
  & \sum_{a} \sum_{\al=1}^{n+1} \langle \psi_\la, e(z_1) \cdots (
  -r_{a;-2\la}(z_\al) h(z_\al) + k dr_{a;-2\la}(z_\al) ) \cdots
  e(z_{n+1}) v_{n} \rangle \omega_{a;2\la}(z) \\ & \nonumber + \sum_{i}
  \sum_{\al=1}^{n+1} \langle \psi_\la, e(z_1) \cdots ( -r_i^{in}(z_\al)
  h(z_\al) + k dr_i^{in}(z_\al) ) \cdots e(z_{n+1}) v_{n} \rangle
  \omega_{i;2\la}^{out}(z).
\end{align}
because of the relation
$$
[f[\eps], e(z)] = -\eps(z) h(z) + k d\eps(z),   
$$ and because we have $f[r_i^{in}]v_n = f[r_{a;2\la}]v_n =0$; the
latter equality is because the $r_{a;2\la}$ have poles of order $\leq
g-1$ at $P_0$.

(\ref{interm}) is then equal to
\begin{align} \label{interm'}
  & \sum_{\al = 1}^{n+1} [-G_{2\la}(z,z_\al)] \langle \psi_\la ,
  h(z_\al) e(z_1)\cdots \check\al \cdots e(z_{n+1}) v_{n}
  \rangle \\ & \nonumber + \sum_{\al = 1}^{n+1} k
  d_{z_\al}G_{2\la}(z,z_\al)] \langle \psi_\la , e(z_1)\cdots
  \check\al \cdots e(z_{n+1}) v_{n} \rangle.
\end{align}
Applying Lemma \ref{lemma:act:h} to the first sum, one gets (\ref{act:f}).
\hfill \qed \medskip

\subsubsection{Action of the Sugawara tensor on the correlation functions}

Let us compute now 
$$
\langle \psi_{\la}, h(z) h(z') e(z_1) \cdots e(z_n) v_n \rangle . 
$$
This is equal to 
\begin{align*}
& \sum_a \langle \psi_{\la}, h[r_a] h(z') e(z_1) \cdots e(z_n) v_n \rangle 
\omega^a(z) 
\\ & + 
\sum_i \langle \psi_{\la}, h[r_{i,0}^{in}] h(z') e(z_1) \cdots e(z_n) v_n 
\rangle \omega_i^{out}(z) 
\end{align*}
that is 
\begin{align*} 
& \sum_a \omega^a(z) \pa_{\la_a}\langle \psi_{\la}, h(z') 
e(z_1) \cdots e(z_n) v_n \rangle 
\\ & + \sum_i \langle \psi_{\la}, h(z')
 h[r^{in}_{i,0}] e(z_1) \cdots e(z_n) v_n \rangle \omega_i^{out}(z) 
\\ & + 2 k d_{z'}G(z,z') \langle \psi_{\la}, e(z_1) 
\cdots e(z_n) v_n \rangle .
\end{align*}
The second line is equal to 
$$
\sum_{\al = 1}^n 2 G(z,z_{\al})
\langle \psi_{\la}, h(z') e(z_1) \cdots e(z_n) v_n \rangle. 
$$
Applying Lemma \ref{lemma:act:h} to the first two sums, we find 
\begin{align*}
& 2(k+2)\langle \psi_{\la}, h(z) h(z') e(z_1) \cdots e(z_n) v_n \rangle 
= 
2 k d_{z'}G(z,z') f_{\la}(z_1,\cdots,z_n) 
\\ & \nonumber + 
\left(\sum_a \omega_a(z) \pa_{\la_a} + 2 \sum_{\al = 1}^n G(z,z_{\al}) 
\right)^2 
f_{\la}(z_1,\cdots,z_n) + O(z,z'). 
\end{align*}

On the other hand, we have, by Lemma \ref{lemma:act:f},
\begin{align*}
  \langle \psi_{\la} , & f(z') e(z) e(z_1) \cdots e(z_n) v_n \rangle
   \rangle \\ & = - G_{2\la}(z',z) (\sum_a
  \omega_a(z)\pa_{\la_a} + \sum_{\al = 1}^n 2G(z,z_\al))
  f_{\la}(z_1,\cdots,z_n) \\ & - \sum_{\al = 1}^n G_{2\la}(z',z_{\al})
  [\sum_a \omega_a(z_{\al})\pa_{\la_a} + \sum_{\beta = 1, \beta\neq
    \al}^n 2G(z_\al,z_\beta)] f_{\la}(z_1,\cdots, z, \cdots, z_n) \\ &
  -\sum_{\al = 1}^n G_{2\la}(z',z_{\al})2G(z_\al,z)
  f_{\la}(z_1,\cdots, z, \cdots, z_n) \\ & + k d_z(G_{2\la}(z',z))
  f_\la(z_1,\cdots,z_n) \\ & + k \sum_{\al = 1}^n
  d_{z_\al}(G_{2\la}(z',z_\al)) f_\la(z_1,\cdots,z,\cdots ,z_n)
\end{align*}
so that
\begin{align*} 
& \langle \psi_{\la}, \left( e(z) f(z') + f(z) e(z')\right)  
e(z_1) \cdots e(z_n) v_n \rangle 
\\ & = 
(\sum_a D_z^{(2\la)}\omega_a(z)\pa_{\la_a} + \sum_{\al = 1}^n 2
( D_z^{(2\la)}\otimes 1)G(z,z_\al))
f_{\la}(z_1,\cdots,z_n)
\\ & 
- 2 \sum_{\al = 1}^n
 G_{2\la}(z,z_{\al}) [\sum_a \omega_a(z_{\al})\pa_{\la_a} + 
\sum_{\beta = 1, \beta\neq \al}^n 2G(z_\al,z_\beta)]
f_{\la}(z_1,\cdots, z, \cdots, z_n)
\\ & 
- 4 \sum_{\al = 1}^n G_{2\la}(z,z_{\al})G(z_\al,z)
f_{\la}(z_1,\cdots, z, \cdots, z_n)
\\ &
+ k [ d_z(G_{2\la}(z',z)) + d_{z'}(G_{2\la}(z,z')) ]
f_\la(z_1,\cdots,z_n)
\\ & 
+ 2 k \sum_{\al = 1}^n d_{z_\al}(G_{2\la}(z,z_\al)) f_\la(z_1,\cdots,z,\cdots
,z_n)
+ O(z - z')
\end{align*}
(with $z$ in the $\al$th place in the r.h.s.), 
where $D^{(\la)}_{z}(\omega)$ is defined by 
\begin{equation} \label{D:la}
D^{(\la)}_{z}(\omega)(z) =  - \lim_{z\to z'}( \omega(z')G_{\la}(z,z') 
+ \omega(z)G_{\la}(z',z) ); 
\end{equation}
$D^{(\la)}_{z}$ defines a connection from the
bundle $\Omega_X$ to $\Omega_X^2$.  

Set 
\begin{align} \label{Tz}
  & (T_z f_\la) (z_1,\ldots,z_n)\\ & = \nonumber \left[ {1\over 2}
    \left(\sum_a \omega_a(z) \pa_{\la_a} + 2 \sum_{\al} G(z,z_{\al})
       \right)^2 \right.  \\ & \nonumber \left.  + \sum_a
      D^{(2\la)}_{z} \omega_a(z) \pa_{\la_a} + 2\sum_{\al}
      (D^{(2\la)}_{z} \otimes 1) (G(z,z_\al)) +
    k \omega_{2\la}(z) \right] f_\la(z_1,\cdots,z_n) \\ & \nonumber +
  \sum_{\al = 1}^n \left[ - 2G_{2\la}(z,z_{\al}) \left( \sum_a
      \omega_a(z_\al) \pa_{\la_a} + 2 \sum_{\beta \neq
        \al}G(z_\al,z_\beta) \right) \right.  \\ & \nonumber
  \left.  + \left(- 4 G_{2\la}(z,z_\al)G(z_\al,z) +2 k
      d_{z_\al}G_{2\la}(z,z_\al) \right) \right]
  f_\la(z_1,\cdots,z,\cdots,z_n), 
\end{align}
where $z$ is in the $\al$th position in the right hand side and we set 
\begin{equation} \label{omega:la}
\omega_{\la}(z) = \limm_{z\to z'} \left( d_{z'} G_{\la}(z,z')
  + d_{z} G_{\la}(z',z) - 2d_{z'} G(z,z') \right) .
\end{equation}

Then 

\begin{prop} \label{thl}
We have 
$$ \langle \psi_\la, T_{\wt \omega}(z) e(z_1)\cdots e(z_n) v_n \rangle
= (T_z f_\la)(z_1,\cdots,z_n).
$$
\end{prop}

\begin{remark}
  It would be interesting to have an expression of the action of
  $T(z)$ directly in terms of the $f_{FS}(z_1,\ldots,z_n)$. For this,
  one would need either to understand the correspondence of Rem.
  \ref{vopros}, or how to express the $T[z^p]v_n$ as combinations of
  the $e[z^{i_1}] \cdots e[z^{i_l}]v_{n+l}$.
\end{remark}

\subsection{Action of the Sugawara tensor in the general case}

In this section, we show how the expression of the operators $T_z$ is
modified in the case of a general semisimple $\bar\G$. For any $\al$
in $\Delta_+$, let $e_\al,f_\al$ and $\al^\vee$ be in $\bar\G_\al$,
$\bar\G_{-\al}$ and $\bar\HH$ forming a standard $\SL_2$-triple, and
let $(a_{ij})_{1\leq i,j\leq r}$ be the Cartan matrix of $\bar\G$.

For $i_1,\ldots,i_s$ in $\{1,\ldots,r\}$, such that $\sum_{j=1}^s
\al_{i_j} < \al$, define the number $n_{\al;i_1,\ldots,i_s}$ by the
equality
$$ [[[f_\al,e_{i_1}],e_{i_2}]\cdots,e_{i_s}] = n_{\al;i_1\ldots i_s}
f_{\al - \sum_{j=1}^s \al_{i_j}};
$$ for $\al,\beta$ in $\Delta_+$, such that $\al-\beta$ belongs to
$\Delta_+$, define the number $N_{\al\beta}$ by the equality 
$$ [f_{\al-\beta},e_\al] = N_{\al\beta} e_{\beta}; 
$$
define $\nu_{i_1\ldots i_s}$ by the equality 
$$ [[e_{i_1},e_{i_2}], \cdots,e_{i_s}] = \nu_{i_1\ldots i_k}
e_{\sum_{j=1}^s e_{\al_{i_j}}}. 
$$

As we have seen, one may attach to a $\G^{out}$-invariant form $\psi$
on any $\G$-module $V$, the forms
\begin{equation} \label{corr:h:rk}
  f_\la(z^{(i)}_u)_{1\leq i\leq r,u\in I_i} = \langle \psi,
  e^{\sum_{a,i} \la^{(i)}_a h_i[r_a]} \prod_{i=1}^r \prod_{u\in I_i}
  e_i(z^{(i)}_u) dz^{(i)}_u v \rangle ,
\end{equation} 
where the $I_i$ are finite sets attached to $i = 1, \ldots,r$ and $v$
is a vector in $V$ with the suitable weight. The $\la^{(i)}_a$ are
formal variables. We attach to them the family $(\la_a)_{1\leq a \leq
  g}$ of formal elements of $\bar\HH^g$, where $\la_a = \sum_i
\la_a^{(i)} h_i$. For $\mu$ in $\bar\HH^*$, we set $(\mu,\la) =
(\mu,\la_a)_{1\leq a \leq g}$.

The form $f_\la(z^{(i)}_u)$ depends on the $z^{(i)}_u$ as a section of
$\Omega_X \otimes \cL_{-(\al_i,\la)}$, regular on $X$ outside $P_0$
and the $z^{(j)}_v$ for the $j$ such that $a_{ij}< 0$. It is
symmetric in the $z^{(i)}_u$ for each $i$, with simple poles at the
diagonals $z^{(i)}_u = z^{(j)}_v$ when $a_{ij}< 0$, and satisfies
$$ \res_{z_{u_1}^{(i)} = z_{v}^{(j)}} \res_{z_{u_2}^{(i)} =
  z_{v}^{(j)}} \cdots \res_{z_{u_{1-a_{ij}}}^{(i)} = z_{v}^{(j)}}
f_\la(z^{(i)}_u) = 0
$$ for $v$ in $I_j$ and $u_1,\ldots,u_{1-a_{ij}}$ distinct in $I_j$
(see \cite{FS}); this is a translation of the Serre relations, using
the identities $\res_{z=z'} \langle \psi, x(z) y(z') v \rangle =
\langle \psi, [x,y](z')v \rangle$.

Assume that $v$ is annihilated by the positive Cartan modes and the
$f_i[z^{i}],i\geq - (g-1)$; let $(h_\nu)_{1\leq \nu \leq r}$ be an
orthonormal basis of $\bar\HH$ and define the Sugawara tensor as
$$ 2(k+h^\vee)T_{\wt \omega}(z) = \limm_{z'\to z}\left( \sum_\nu h_\nu(z)
  h_\nu(z') + \sum_{\al\in\Delta_+} f_\al(z) e_\al(z') + e_\al(z)
  f_\al(z') - k (\dimm\bar\G) \wt\omega(z,z')\right),  
$$ with $h^\vee$ the dual Coxeter number of $\bar\G$.

Let $P$ (resp. $P'$) be the set of sequences $p = (i_1,\ldots ,i_s)$
such that $\al = \sum_{j=1}^s \al_{i_j}$ (resp.  $\al>\sum_{j=1}^s
\al_{i_j}$). The sequence $(u_j)$ is associated to $P$ if it is a
sequence of pairwise different elements of $\cup_i I_i$, such that
$u_k$ belongs to $I_{i_k}$. We denote by $S_i$ the subset of $I_i$
formed by all $u_j$ such that $i_j$ is equal to $i$.

\begin{prop}
The action of $T_{\wt\omega}(z)$ on the correlation function
(\ref{corr:h:rk}) is given by 
\begin{align} \label{Tz:hrk}
  & 2(k+h^\vee) \langle \psi_\la, T_{\wt \omega}(z) \prod_{i=1}^r
  \prod_{u\in I_i} e_i(z^{(i)}_u) v \rangle \\ & \nonumber = \left[
    \sum_\nu \left( \sum_{a} \omega_a(z) \pa_{(h_\nu)_a} +
      \sum_{i}\sum_{u\in I_i} (\al_i,h_\nu) G(z,z^{(i)}_{u}) \right)^2
  \right. \\ & \nonumber \left. + \sum_{\al\in \Delta_+} \left( \sum_a
      D_z^{(\al,\la)} \omega_a(z) \pa_{(\al^\vee)_a} + \sum_{i}
      \sum_{u\in I_i}(\al_i,\al^{\vee}) D_z^{(\al,\la)} G(z,z^{(i)}_u)
    \right) + k \sum_{\al\in \Delta_+} \omega_{(\al,\la)}(z) \right]
  \\ & \nonumber f_{\la}(z^{(i)}_u) \\ & \nonumber +
  \sum_{p'\in P'} 
  n_{\al;i_1\ldots i_s} N_{\al;\al_{i_1} + \cdots + \al_{i_s}} /
  \nu_{i_1 \ldots i_s} \\ & \nonumber \sum_{(u_i)\ \on{associated}\ 
    \on{to}\ p'} G_{(\al,\la)}(z,z^{(i_1)}_{u_1})
  G_{(\al-\al_{i_1},\la)}(z^{(i_1)}_{u_1},z^{(i_2)}_{u_2}) \cdots
  G_{(\al-(\al_{i_1}+ \cdots + \al_{i_s}),\la)}(z^{(i_s)}_{u_s},z) \\ 
  & \nonumber \res_{z^{(i_1)}_{u_1} =
    z^{(i_2)}_{u_2}}\res_{z^{(i_2)}_{u_2} = z^{(i_3)}_{u_3}} \cdots
  \res_{z^{(i_{s-1})}_{u_{s-1}} = z} f(z^{(i)}_u)_{|z^{(i_s)}_{u_s} =
    z} \\ & \nonumber -
\sum_{p\in P}  
n_{\al;i_1 \ldots i_{s-1}} / \nu_{i_1 \ldots i_s} \\ & \nonumber
\sum_{(u_i)\ \on{associated}\ \on{to} \ p}
G_{(\al,\la)}(z,z^{(i_1)}_{u_1}) G_{(\al -
  \al_{i_1},\la)}(z^{(i_1)}_{u_1}, z^{(i_2)}_{u_2}) \cdots G_{(\al -
  (\al_{i_1}+ \cdots + \al_{i_{s-1}}),\la)}(z^{(i_{s-1})}_{i_{s-1}},
z^{(i_s)}_{u_s}) \\ & \nonumber [\sum_a \omega_a(z^{(i_s)}_{u_s})
\pa_{(\al_{i_s}^\vee)_a} + \sum_{i} \sum_{u \in I_i - S_i}
(\al_i,\al_{i_s}^\vee) G(z^{(i_s)}_{u_s}, z^{(i)}_{u}) +
(\al,\al_{i_s}^\vee)G(z^{(i_s)}_{u_s},z) ] \\ & \nonumber
\res_{z^{(i_1)}_{u_1}= z^{(i_2)}_{u_2}} \res_{z^{(i_2)}_{u_2}=
  z^{(i_3)}_{u_3}} \ldots \res_{z^{(i_{s-1})}_{u_{s-1}}= z}
f_\la(z^{(i)}_u)_{|z^{(i_s)}_{u_s} = z} \\ & \nonumber +
\sum_{p\in P}  
n_{\al;i_1 \ldots i_{s-1}} / \nu_{i_1 \ldots i_s} \\ & \nonumber
\sum_{(u_i)\ \on{associated} \ \on{to} \ p}
G_{(\al,\la)}(z,z^{(i_1)}_{u_1}) G_{(\al -
  \al_{i_1},\la)}(z^{(i_1)}_{u_1}, z^{(i_2)}_{u_2}) \cdots G_{(\al -
  (\al_{i_1}+ \cdots + \al_{i_{s-1}}),\la)}(z^{(i_{s-1})}_{u_{s-1}},
z^{(i_s)}_{u_s}) \\ & \nonumber k d_{z} G_{(\al - (\al_{i_1} + \cdots
  + \al_{i_s}),\la)} (z^{(i_{s-1})}_{u_{s-1}}, z)
\res_{z^{(i_1)}_{u_1}= z^{(i_2)}_{u_2}} \res_{z^{(i_2)}_{u_2}=
  z^{(i_3)}_{u_3}} \ldots \res_{z^{(i_{s-1})}_{u_{s-1}}= z}
f_\la(z^{(i)}_u)_{|z^{(i_s)}_{u_s} = z} ,
\end{align}
where for $x$ in $\bar\HH$, we denote by $x_a$ the element
$(0,\ldots,x,\ldots,0)$ of $\bar\HH^g$ ($x$ at the $a$th place); and
by $\pa_h$ the partial derivative in $\bar\HH^g$ in the direction of
$h$, for $h$ in $\HH^g$.
\end{prop}

\begin{remark}
  The set $P$ and its associated sequences appeared in the work
  \cite{Sch-Varch} on integral formulas for the KZ equations.
\end{remark}

\subsection{Expression of the KZB connection} \label{KZB}

Denote by $Proj_g^{(1)}$ the moduli space of quadruples $\wt m =
(X,[\{\z_\al\}],P_0,z)$, where $X$ is a curve of genus $g$,
$[\{\z_\al\}]$ is a projective atlas of $X$ (that is an atlas whose
transition functions are projective transformations), $P_0$ a point of
$X$ and $z$ a coordinate of the atlas with origin at $P_0$. A local
coordinate related to some $z_\al$ by a projective transformations
will be called a projective coordinate.

For each representation $V$ of $\G^{out}$, we may form the bundle
$CB(V)$ over $Proj_g^{(1)}$, whose fiber at $\wt m$ is defined as the
space of $\G^{out}$-invariant forms on $V$.

A projectively flat connection on the bundle $CB(V)$ is defined as
follows. Let $\wt m\mapsto \psi(\wt m)$ be a local section of $CB(V)$.
Let $\delta\wt m$ be a variation of $\wt m$. Then
\begin{equation} \label{connect}
  \nabla_{\delta\wt m}\psi = \pa_{\delta\wt m} \psi - \psi \circ
  T_0[\xi(\delta\wt m)],
\end{equation}
where the equality is in $V^*$ and $\xi(\delta\wt m)$ is the element
of $\CC((z))\pa_z$ induced by $\delta\wt m$ (for any moduli $\wt m$,
we have a ring $R_{\wt m}$ contained in $\CC((z))$, and we set $R_{\wt m
  + \delta\wt m} = (1 + \xi(\delta\wt m))R_{\wt m}$). We set $T_0[\xi]
= \res_{P_0}(T_0(z) dz^2 \xi(\delta \wt m)(z) \pa_z)$, with $T_0(z)$
defined as $T_{\wt \omega}(z)$ in (\ref{expr:T}) replacing $\wt
\omega$ by ${{dzdw}\over{(z-w)^2}}$.

This connection is well-defined, preserves $CB(V)$ and is projectively
flat (see \cite{TUY}).

The form $\wt \omega$ defined by (\ref{wt:omega}) depends only on the
choice of $a$-cycles. On the other hand, this form determines a
projective structure on $X$. Indeed, it is known that there is a
bijective correspondence between bidifferential forms near the
diagonal with behavior ${{dzdw}\over{(z-w)^2}} + r(z) dzdw +
o(z-w)dzdw$, up to addition of regular bidifferential forms vanishing
on the diagonal, and projective structures on $X$. The correspondence
associates to the projective atlas $[\{\z_{\al}\}]$ the form
$d_{\z_\al} d_{\z'_\al} \ln(\z_\al - \z'_\al)$. Conversely, the
projective coordinate $\zeta$ associated to the bidifferential form
${{dzdz'}\over{(z-z')^2}} + r(z) dzdz' + o(z-z')dzdz'$ is determined
by the equation $S(\zeta,z) = -6r(z)$, where $S(\zeta,z)$ is the
Schwarzian derivative of $\zeta$ with respect to $z$. Then
$T_0(\zeta)(d\zeta)^2$, computed in a projective coordinate determined
by $\wt \omega$, gets identified with $T_{\wt\omega}(z)(dz)^2$.

Let us define $\cM_g^{(a)}$ as the moduli space of genus $g$ curves
with marked homology classes of $a$-cycles. $\wt \omega$ defines a map
from $\cM_g^{(a)}$ to $Proj_g$, such that its composition with
projection of $Proj_g$ to $\cM_g$ coincides with the projection of
$\cM_g^{(a)}$ on $\cM_g$.

Define $\cM_g^{(a)(1)}$ as the fibered product of $\cM_g^{(a)}$ with
$Proj_g^{(1)}$ over $Proj_g$. The KZB connection is defined on
$Proj_g^{(1)}$, and it induces a connection on $\cM_g^{(a)(1)}$, using
the map from $\cM_g^{(a)(1)}$ to $Proj_g^{(1)}$. This connection can
be expressed as follows.

Let us express the connection induced by (\ref{connect}) in terms of
correlation functions.  For any formal vector field $\xi = \xi(z)
\pa_z$ in $\CC((z))\pa_z$, set $\eps^2 = 0$ and $R_\eps =
(1+\eps\xi)R$; let $\Omega_R \subset \Omega_{\cK}$ be the space of
differentials of $R$ and $\Omega_{R_\eps}$ the space of differentials
of $R_\eps$.  Then $\Omega_{R_\eps}$ is equal to
$(1+\eps\cL_\xi)(\Omega_R)$, where $\cL_\xi$ is the Lie derivative
associated to $\xi$. Similarly, we have $dR_\eps =
(1+\eps\cL_\xi)(dR)$. Therefore, $1+\eps\cL_\xi$ induces a map from
$\Omega_R/dR$ to $\Omega_{R_\eps}/dR_\eps$. Bases of these spaces are
the classes of the $\omega_a$ and $dr_a$. On the other hand, we have
the formula $\int_{\gamma'}(1+\eps\cL_\xi)(\omega) = \int_\gamma
\omega$ for any cycle $\gamma$ of $X$, deformed to $\gamma'$ and any
$\omega$ in $\Omega_R$. Therefore, we have $\cL_\xi (dr_a) = 0$ mod.
$dR$ and $\cL_{\xi} \omega_a = \sum_b \delta\tau_{ab} dr_b$ mod $dR$,
where $\delta\tau_{ab}$ is the variation of the period matrix.

{}From this follows: 

\begin{prop}
  Let $\wt m \mapsto \psi(\wt m)$ be a section of the bundle
  $\cF^{(n)}(m')$ over $\cM_g^{(a)(1)}$, then the KZB connection is
  expressed as
\begin{align*}
  & \nabla_{\delta\wt m} f(\wt m)_\la(z_1,\ldots,z_n) = \pa_{\delta\wt
    m} f(\wt m)_\la(z_1,\ldots,z_n) - \langle \psi_\la, T[\xi(\delta
  \wt m)] e(z_1) \cdots e(z_n) v_n \rangle ,
\end{align*}
where $\langle \psi_\la, {1\over{k+2}} T[\xi(\delta \wt m)] e(z_1)
\cdots e(z_n) v_n \rangle$ can be computed using (\ref{Tz}).
\end{prop}

\begin{remark} 
  The fact that the action of $T(z)$ preserves the vanishing
  conditions of Feigin-Stoyanovsky (vanishing on codimension $k$
  diagonals) probably again follows from the identity $(e^k)' =
  :he^k:$.
\end{remark}

\subsection{Motion of marked points ($\SL_2$ case)} \label{marked}

In this section, we indicate how the above results are changed in the
case of marked points. Let $(P_i)_{i = 1 ,\ldots,N}$ be marked points
on $X$, distinct from $P_0$. Attach to each $P_i$ the weight $\La_i$
and the evaluation Verma module $V_{\La_i}$. $V_{\La_i}$ is generated
by the vector $v_{-\La_i}$ such that $h v_{-\La_i} = -\La_i v_{-\La_i}$,
and $f v_{-\La_i} = 0$. Set again $\psi_\la = \psi \circ e^{\sum_a
  \la_a h[r_a]}$ and
$$ f_\la (z_1,\ldots,z_m) = \langle \psi_\la, \left( e(z_1)dz_1\cdots
  e(z_m)dz_m v_n\right) \otimes v_{-\La_1} \otimes \cdots \otimes v_{-\La_N}
\rangle,
$$ $m = n -{1\over 2} \sum_i \La_i$.

$f_\la(z_1,\cdots,z_m)$ depends on the $z_\al$ as a section of
$\Omega_X\cL_{2\la}$, regular outside $P_0$ and with simple poles at
the $P_i$.

For $w_i$ in $V_{\La_i}$, the values of the $\langle \psi_\la,
(e(z_1)\cdots e(z_m)v_{n}) \otimes (\otimes_{i=1}^N w_i) \rangle$ can
be recovered from $f_\la(z_1,\cdots,z_m)$ using the rule
\begin{align*}
  & \res_{z = P_i} \langle \psi_\la , (e(z)dze(z_1)dz_1\cdots
  e(z_m)dz_m v_{n}) \otimes (\otimes_{i=1}^N w_i) \rangle \\ & = -
  \langle \psi_\la , (e(z_1)dz_1 \cdots e(z_m)dz_m v_{n}) \otimes
  e^{(i)}(\otimes_{i=1}^N w_i) \rangle.
\end{align*}

The action of the Sugawara tensor is expressed as 
\begin{align*}
  & \langle \psi_{\la} , \left( T_{\wt\omega}(z) e(z_1)dz_1 \cdots
    e(z_m)dz_m v_n \right) \otimes (\otimes_{i=1}^N v_{-\La_i}) \rangle
  \\ & = \left[ {1\over 2} \left(\sum_a \omega_a(z) \pa_{\la_a} + 2
      \sum_{\al} G(z,z_{\al}) - \sum_i \La_i G(z,P_i)\right)^2
  \right.\\ & \left.  + \sum_a D^{(2\la)}_{z} \omega_a(z) \pa_{\la_a}
    + 2\sum_{\al} (D^{(2\la)}_{z} \otimes 1) G(z,z_\al) - \sum_{i}
    \La_i (D^{(2\la)}_{z} \otimes 1) G(z,P_i) + k \omega_{2\la}(z)
  \right] \\ & f_\la(z_1,\cdots,z_n) \\ & + \sum_{\al = 1}^n \left[ -
    2G_{2\la}(z,z_{\al}) \left( \sum_a \omega_a(z_\al) \pa_{\la_a} + 2
      \sum_{\beta \neq \al}G(z_\al,z_\beta) - \sum_{i}\La_i
      G(z_\al,P_i) \right) \right.  \\ & \left.  + \left(- 4
      G_{2\la}(z,z_\al)G(z_\al,z) +2 k d_{z_\al}G_{2\la}(z,z_\al)
    \right) \right] f_\la(z_1,\cdots,z,\cdots,z_n) .
\end{align*}

When $k = -2$, the r.h.s. of this formula is the expression for a
commuting family of differential-difference operators, or
alternatively, for a commuting family of differential operators acting
on some finite-dimensional bundle over $J^0(X)$.

The KZB connection is now a connection over the bundle of conformal
blocks over $Proj_g^{(n)}$, which is the set of quadruples $\wt m =
(X,[\{\z_\al\}],P_i,\z_i)$ of curves with projective structure, $n$
marked points and flat coordinates vanishing at these points.

The vector fields $\z_i {{\pa}\over{\pa \z_i}}$ describing the changes
of coordinates fixing the points, and ${{\pa}\over{\pa P_i}}$
describing the changes of points in the fixed coordinate, are
respectively given by the action of Sugawara elements corresponding to
vector fields $\xi_{{{\pa}\over{\pa P_i}}}$ equal to ${{\pa}\over{\pa
    \z_i}}$ at $P_i$ and $o(\z_j)$ at $P_j$, and $\xi_{\z_i
  {{\pa}\over{\pa \z_i}}}$ equal to $\z_i {{\pa}\over{\pa \z_i}}$ at
$P_i$ and $o(\z_j)$ at $P_j$.

Set $G(z,w) = {{dz}\over {z-w}} + \phi(z) dz + o(z-w)dz$, so that
$G_\la(z,w) dz = {dz\over{z-w}} + g_\la(z) dz$, with 
\begin{align*} & g_\la(z)dz = \\ & \phi(z) dz + \sum_{a=1}^g \omega_a(z) \left(
    \pa_{\eps_a}\ln\Theta(-\la+(g-1)A(P_0)-\Delta) -
    \pa_{\eps_a}\ln\Theta(gA(P_0)-A(z)-\Delta) \right).
\end{align*}
\begin{prop}
  The KZB connection is expressed, in the direction of variation of
  coordinates at $P_i$ by
  $$ 2(k+2) \nabla_{\z_i {{\pa}\over{\pa \z_i}}} f(\wt m)_{\la}(z_\al)
  = 2(k+2) \z_i {{\pa}\over{\pa \z_i}}f(\wt m)_{\la}(z_\al) - {1\over
    2} \La_i(\La_i +2)f(\wt m)_{\la}(z_\al).
$$ and in the direction of variation of $P_i$, by
\begin{align} \label{var:pts}
  & 2(k+2)\nabla_{{{\pa}\over{\pa P_i}}} f(\wt m)_\la(z_\al) = 2(k+2)
  {{\pa}\over{\pa P_i}} f(\wt m)_\la(z_\al) - [ - \La_i\sum_a
  \omega_a(P_i)\pa_{\la_a} \\ & \nonumber + \La_i \left( \sum_{j\neq
      i} \La_j G(P_i,P_j) - 2 \sum_\al G(P_i,z_\al) \right) + \La_i^2
  \phi(P_i) + 2\La_i g_{2\la}(P_i) ] f(\wt m)_\la \\ & \nonumber +
  \sum_{\al}[-2 G_{2\la}(P_i,z_\al)(\sum_a \omega_a(z_\al) \pa_{\la_a}
  + 2 \sum_{\beta \neq \al} G(z_{\al},z_\beta)) \\ & \nonumber -4
  G_{2\la}(P_i,z_\al)G(z_\al,P_i) + 2 k d_{z_\al} G_{2\la} (z_\al,P_i)
  ] \res_{z = P_i} f(\wt m)_\la (z,(z_\beta)_{\beta\neq \al}) ;
\end{align}
when $m = 0$, this equation simplifies to   
\begin{align} \label{var:pts:0}
  & 2(k+2)\nabla_{{{\pa}\over{\pa P_i}}} f(\wt m)_\la = 2(k+2)
  {{\pa}\over{\pa P_i}} f(\wt m)_\la \\ & \nonumber - [ - \La_i\sum_a
  \omega_a(P_i)\pa_{\la_a} + \La_i\sum_{j\neq i} \La_j G(P_i,P_j) +
  \La_i^2 \phi(P_i) + 2\La_i g_{2\la}(P_i) ] f(\wt m)_\la.
\end{align}
\end{prop}

\begin{remark}
  It would be interesting to express the equations obtained above in
  terms of dynamical $r$-matrices, as it was done in \cite{Houches}.
\end{remark}

\section{Commuting differential operators} \label{comm:do} 

It is natural to consider the operators (\ref{Tz}) as
differential-evaluation operators acting on functions on $J^0(X)^r
\times \prod_{i=1}^r S^{n_i}X$, for $k$ an arbitrary number. In
particular, when $k$ is critical, one expects these operators to
commute with each other. To prove this, we will consider modules
$W_{n|m,m'}$ generalizing the twisted Weyl modules.

For generic $\la_0$ in $J^0(X)$, $\la_0$-twisted conformal blocks for
these modules can be characterized via functions (\ref{oo}) as
formal sections of finite-dimensional bundles over $J^0(X)$.

\subsection{Twisted conformal blocks for general modules} \label{formal}

Let $X$ be a smooth complex curve of genus $g\geq 1$ and let $P_0$ be
a fixed point of $X$. Denote by $\cK$ and $\cO$ the local field and
ring of $X$ at $P_0$. Denote also by $R$ the ring $H^0(X -
\{P_0\},\cO_X)$ and by $\AAA$ the adeles ring of $X$.


Recall that (\ref{formal:psi:la}) defined a form $\psi_\la$, depending
on formal variables $\la^{(i)}_a$, on an arbitrary $\G$-module $V$.

For $\mu_1,\ldots,\mu_g$ complex linear combinations of the
$\la_a^{(i)}$, define $R_{(\mu_i)}^{(f)}$ as the subspace of
$\cK[[\la_a^{(i)}]]$ formed by the functions $f(z,\la_a^{(i)})$
depending formally on the $\la_a^{(i)}$, such that the coefficients of
the monomials in $\la_a^{(i)}$ extend to regular functions on $\wt X -
\si^{-1}(P_0)$ and we have $f(\gamma_{A_a}z,\la_a^{(i)}) =
f(z,\la_a^{(i)})$ and $f(\gamma_{B_a}z,\la_a^{(i)}) =
e^{\mu_a}f(z,\la_a^{(i)})$.  

$\psi_\la$ has the following properties:

\begin{lemma}
  a) Set for $a = 1,\ldots,g$, $\la_a = \sum_{i} \la_a^{(i)} h_i$.
  Define $\G_\la^{out(f)}$ as
  $$ \G_\la^{out(f)} = (\bar\HH \otimes R)[[\la_a^{(i)}]] \oplus
  \oplus_{\al\in \Delta} (\bar\G_\la \otimes R^{(f)}_{\langle \al,\la_1
    \rangle , \ldots , \langle \al,\la_g \rangle }).
  $$ Then $\psi_\la$ is $\G^{out(f)}_\la$-invariant.

  b) $\la\mapsto \langle \psi_\la,v\rangle$ satisfies the differential
  equation $\pa_{\la^{(i)}_a}\langle \psi_\la , v \rangle = \langle
  \psi_\la, h_i[r_a] v \rangle$ for any $v$ in $V$.
\end{lemma}

{\em Proof.} Clearly, $\G^{out(f)}_\la$ is contained in
$\Ad(e^{-\sum_{i,a}\la^{(i)}_a
  h_i[r_a]})(\G^{out}[[\la_a^{(i)}]][\la_a^{(i)-1}])$; this implies
a).  We have for any $a,b = 1, \ldots,g$, $\langle dr_a, r_b\rangle =
{1\over{2i\pi}}\int_{\pa i(X)} dr_a r_b$; the contributions of the
paths $\wt B_c$ and $\wt B_c^{-1}$ cancel each other, as well as those
of the paths $\wt A_c$ and $\wt A_c^{-1}$, $c\neq b$; the sum of the
contributions of the paths $\wt A_b$ and $\wt A_b^{-1}$ is equal to
${1\over{2i\pi}}\int_{A_b} dr_a$, which is zero as $r_a$ is
single-valued along $a$-cycles. Therefore we have $[h_i[r_a],h_j[r_b]]
= 0$ for any $i,j,a,b$, which proves b).  \hfill \qed \medskip

\subsection{Conformal blocks for the $W_{n|m,m'}$}
\label{sect:tw:corr}

In this section, we set $\bar\G = \SL_2$.  Let $k$ be an arbitrary
complex number.

For $m,m'$ integer numbers with $m+m' \ge 0$, define $\G^{in}_{m,m'}$ by 
$$\G^{in}_{m,m'} = (\bar\N_- \otimes z^{m}\cO) \oplus (\bar\HH \otimes
\cO) \oplus (\bar\N_+ \otimes z^{m'}\cO) \oplus \CC K.
$$ Define $\G^{in}_{m,\infty}$ and $\G^{in}_{-\infty,\infty}$ by the
convention that $z^\infty \cO = 0$ and $z^{-\infty}\cO = \cK$.

Let $n$ be a positive integer. If $m+m' >0$, $(m,m') =
(-\infty,\infty)$, or $m+m' = 0$ and $n = -km$, define $\chi_{n |
  m,m'}$ as the character of $\G^{in}_{m,m'}$ such that $\chi_{n |
  m,m'}(K) = k$, $\chi_{n | m,m'}(h[z^i]) = -2n \delta_{i,0}k$,
$\chi_{n | m,m'}(x[z^i]) = 0$, $x = e,f$.

Define $W_{n|m,m'}$ as the induced module $U\G
\otimes_{U\G^{in}_{m,m'} } \CC_{\chi_{n|m,m'}}$. Denote by $v_n$ the
vector $1\otimes 1$ of this module. (When $m+m' = 0$, $W_{n|m,m'}$ is
a twisted Weyl module.)  For $\la_0$ a complex number, define
$CB_{\la_0}(W_{n|m,m'})$ as the space of $\G_{\la_0}^{out}$-invariant
linear forms on $W_{n|m,m'}$ (where $\G^{out}_{\la_0}$ is as in Thm.\
\ref{seerose}).

Let us define $\cF_{\la_0}^{(n)}$ as the space of forms
$f_{\la}(z_1,\ldots,z_n)$, depending formally on $\la$
in the neighborhood of $\la_0$, symmetric in $z_1,\ldots, z_n$,
sections of $\Omega_X\cL_{-2\la}$ in $z_i$, regular outside $P_0$.
Define for any integer $p$, $\cF_{\la_0}^{(n)}(p)$ as the subspace of
$\cF_{\la_0}^{(n)}$ consisting of the forms with poles at $z_i = P_0$
of order at most $p$.

For any $\rho$ in $R_{2\la}$, define first order differential
operators $\wt f[\rho]$ by 
\begin{align} \label{wtf}
  & (\wt f[\rho] f_\la)(z_1,\ldots,z_{n+1}) \\ & \nonumber =
  \sum_{i=1}^{n+1} [-\rho(z_i) (\sum_a \omega_a(z_i) \pa_{\la_a} + 2
  \sum_{j\neq i} G(z_i,z_j)) + k d\rho(z_i)] f_\la(z_1,\ldots
  \check i \ldots z_{n+1}).
\end{align}
$\wt f[\rho]$ maps $\cF_{\la_0}^{(n)}$ to $\cF_{\la_0}^{(n+1)}$.

\begin{prop} \label{cb}
  Define a map $\iota$ from $CB_{\la_0}(W_{n|m,m'}) \to
  \cF_{\la_0}^{(n)}$ by 
  $$\iota(\psi_{\la_0})(\la|z_1,\ldots,z_n) = \langle \psi_{\la_0} ,
  e^{\sum_a (\la-\la_0)_a h[r_a] }e(z_1)\cdots e(z_n) v_n \rangle
   , $$ for $\psi_{\la_0}$ in $CB_{\la_0}(W_{n|m,m'})$.
  
  Assume that $H^1(X,\cL_{2\la_0}(-mP_0))$ is zero. Then $\iota$ is an
  isomorphism from $CB_{\la_0}(W_{n| m,m'})$ to the intersection of
  the kernels of the $\wt f[\rho]$ in $\cF_{\la_0}^{(n)}(m')$, with
  $\rho$ in $R_{2\la}\cap z^m \cO$ (which we may view as $H^0(X,
  \cL_{2\la}(-m P_0))$).
\end{prop}

{\em Proof.} That the image of $\iota$ is contained in the kernel of
the $\wt f[\rho]$ follows from the identity
$$ \langle \psi_{\la_0} , e^{\sum_a (\la - \la_0)_a h[r_a]} [f[\rho],
e(z_1)\cdots e(z_n)]v_n \rangle = 0 ,
$$ which follows from $f[\rho]v_n = 0$ and $\langle \psi_{\la_0},
f[e^{-2\sum_a(\la - \la_0)_a r_a}\rho]v \rangle = 0$ for any vector $v$.

Let us now take $f_\la(z_1,\ldots,z_n)$ in $\cF_{\la_0}^{(n)}(m')$, in
the kernel of the $\wt f[\rho]$ and let us construct its preimage by
$\iota$. 

Clearly, $CB_{\la_0}(W_{n|m,m'})$ is isomorphic to the space of linear
forms $\phi$ on $U\G$, such that $\phi(xx^{in}) = \phi(x^{out}x) = 0$,
for $x^{in}$ in $\G^{in}_{m,m'}$ and $x^{out}$ in $\G_{\la_0}^{out}$.

Define $\CC\langle h[r_a],e[\eps] \rangle$ as the subalgebra of $U\G$
generated by the $h[r_a]$ and the $e[\eps]$, $\eps$ in $\cK$.  Since
we have $\cK = R_{2\la_0} + z^m \cO$, the map
$$
\pi : U\G_{\la_0}^{out} \otimes \CC\langle h[r_a],e[\eps] \rangle \otimes 
U\G^{in}_{m,m'} \to U\G
$$ given by the product is surjective. It kernel is spanned by the $a
e[\eps]\otimes b \otimes c - a\otimes e[\eps]b \otimes b$, $\eps$ in
$R_{-2\la_0}$, the $a\otimes b e[\eps] \otimes c - a \otimes b \otimes
e[\eps]c$, $\eps$ in $z^{m'}\cO$, the $ah[1] \otimes b \otimes c - a
\otimes b \otimes h[1]c - a \otimes [h[1],b]\otimes c$ and the
$af[\eps] \otimes b \otimes c - a \otimes b \otimes f[\eps]c - \sum a
[f[\eps],b]' \otimes [f[\eps],b]'' \otimes [f[\eps],b]''' \otimes c$,
$\eps$ in $R_{2\la_0}\cap z^m \cO$ with $a,b,c$ in
$U\G_{\la_0}^{out}$, $\CC\langle h[r_a],e[\eps]\rangle$ and
$U\G^{in}_{m,m'}$, and $\sum [f[\eps],b]'\otimes [f[\eps],b]'' \otimes
[f[\eps],b]'''$ any preimage of $[f[\eps],b]$ by $\pi$.

Define a linear form $\bar\phi$ on $\CC\langle h[r_a],e[\eps] \rangle$
by the formula 
$$\bar\phi(\prod_a h[r_a]^{\al_a} e[\eps_1]\cdots e[\eps_{n'}]) =
\delta_{nn'} \res_{z_1 = P_0}\cdots \res_{z_n = P_0}
f_{(\al_a)}(z_1,\ldots,z_n)\eps_1(z_1)\cdots \eps_n(z_n),
$$ where we set $f_{\la}(z_1,\ldots,z_n) = \sum_{(\al_i)} \prod_a (\la
- \la_0)_a^{\al_a}f_{(\al_a)}(z_1,\ldots,z_n)$.  Extend $\bar\phi$ to
$U\G_{\la_0}^{out} \otimes \CC\langle h[r_a],e[\eps] \rangle \otimes
U\G^{in}_{m,m'}$ by the rule that $\bar\phi(a\otimes b \otimes c) =
\eps(a)\eps(c)\bar\phi(b)$, $\eps$ denoting the counit.

The functional properties of $f_{\la}(z_1,\ldots,z_n)$ imply that the
image of the kernel of $\pi$ is mapped to $0$ by $\bar\phi$, so that
$\bar\phi$ defines a linear form of $U\G$. It is then clear that this
form is left $\G^{out}_{\la_0}$-invariant and right
$\G^{in}_{m,m'}$-invariant, and that its image by $\iota$ is
$f_{\la}(z_1,\ldots,z_n)$.  \hfill \qed\medskip

\begin{lemma} \label{tauto} The operator $T(z)$ acts naturally on 
  $CB_{\la_0}(W_{n|m,m'})$. When $m$ is $\leq - (g-1)$, this action is
  expressed on the $f_\la(z_1,\ldots,z_n)$ by formula (\ref{Tz}).
\end{lemma}

\begin{remark} Since 
  $H^1(X,\cL_{2\la_0}(-mP_0))$ is zero, $H^1(X,\cL_{2\la}(-mP_0))$
  also vanishes for $\la$ in a neighborhood to $\la_0$. By the
  Riemann-Roch theorem, it follows that $H^0(X, \cL_{2\la}(-m P_0))$
  has constant dimension at the neighborhood of $\la_0$.  It follows
  that the $\rho$ understood in the statement of Prop.\ \ref{cb} form a
  free $\CC[[(\la-\la_0)_a]]$-module with rank equal to this
  dimension. \hfill \qed\medskip
\end{remark}

\begin{remark} \label{8} 
  The condition that $H^1(X,\cL_{2\la_0}(-mP_0))$ vanishes is
  fulfilled if $m<-(g-1)$ and any $\la_0$, or if $m = -(g-1)$ and
  $2\la_0$ not one some translate of the theta-characteristic
  containing zero. In the latter case, $CB_{\la_0}(W_{n|m,m'})$ is
  isomorphic to $\cF_{\la_0}^{(n)}(m')$, because
  $H^0(X,\cL_{2\la_0}(-mP_0))$ also vanishes. \hfill \qed\medskip
\end{remark}

\begin{remark} If $m = -(g-1)$ and
  $2\la_0$ is in the translate of the theta-characteristic (for
  example, if $\la_0$ is zero), the image of $\iota$ is characterized
  by some vanishing conditions near $\la_0$.
\end{remark}

\subsection{Commuting differential operators} 

\begin{thm} \label{main}
  Suppose that $k$ equals $-2$.

  1) Set for $p$ integer $\ge g$ and $\la$ in $J^0(X)$,
  $\cF^{(n)}_\la(p) = S^n H^0(X, \Omega_X \cL_{-2\la}(pP_0))$.
  $(\cF^{(n)}_\la(p))_{\la\in J^0(X)}$ forms a finite-dimensional
  vector bundle over $J^0(X)$, denoted $\cF^{(n)}(p)$.  The operators
  $T_z$ defined by (\ref{Tz}) form a family of commuting differential
  operators acting on sections of this bundle.  This family has rank
  $\leq 3g-3 + p$.  It normalizes the first order operators $\wt
  f[\rho]$ defined by (\ref{wtf}), $\rho$ in
  $H^0(X,\cL_{2\la}(-mP_0))$ for any $m$ (that is, it preserves the
  intersection of their kernels). 

  2) Formula (\ref{Tz}) also defines a family of commuting
  differential-evaluation operators, acting on functions of $\la$ in
  $J^0(X)$ and of $z_1,\ldots,z_n$ in a subset $U$ of $X$ (e.g. the
  pointed formal disc at $P_0$), symmetric in the $z_i$; these
  operators are indexed by points of $U$. They normalize the operators
  $\wt f[\rho]$, $\rho$ some function on $U$, defined by formula
  (\ref{wtf}).
\end{thm}

{\em Proof.} Let us prove 1). If $p\ge g$, the action of $T_z$ on the
jets at $\la_0$ of sections of $\cF^{(n)}(p)$ coincides with the
action of $T(z)$ on $CB_{\la_0}(W_{n|-(g-1),p})$, by Rem.  \ref{8} and
Lemma \ref{tauto}.  Since the $T(z)$ commute together, this shows that
the operators $T_z$ form a commutative family. The result on
normalization of the $\wt f[\rho]$ follows from the fact that the
action of $T_z$ on the intersection of their kernels coincides with
the action of $T(z)$ on $CB_{\la_0}(W_{n|m,\infty})$.

\begin{lemma} For any $f_\la$, $(T_zf_\la)(z_1,\ldots,z_n)$ is a
  quadratic form on $z$, regular on $X$ except for a pole of order
  $\leq p$ at $P_0$.
\end{lemma}

{\em Proof of Lemma.} It is clear that the r.h.s. of (\ref{Tz}) is a
quadratic form in $z$ with possible poles at $P_0$ and the $z_\al$.
Since $k = -2$,  one checks that this expression has no
pole at $z_\al$. 

Let us evaluate the pole at $P_0$. Let $z$ be a local coordinate at
$P_0$. $G_\la(z,w)$ has the expansion 
$$ G_{2\la}(z,w) = {{z^{g-1}w^{1-g}dz}\over{z-w}} + z^{g-1}w^{1-g}dz
\sum_{i,j\ge 0} a_{ij}(\la) z^i w^j. 
$$

Therefore, if $\omega$ belongs to $H^0(X,\Omega_X)$, then
$D_z^{(2\la)}\omega$ is in $H^0(X,\Omega_X^2(P_0))$, because if
$\omega_a$ is $(z^a + o(z^a)) dz$, we have $D_z^{(2\la)} \omega_a =
[(2g-2-a)z^{a-1}+O(z^a))dz$. 

On the other hand, $\omega_{2\la}$ has the expansion at $P_0$ 
$$ \omega_{2\la} = -g(g-1)z^{-2}(dz)^2 - 2 (g-1)z^{-1} (dz)^2 a_{00}(\la)
+ O(1) (dz)^2.
$$
So the two first lines of the r.h.s.  of
(\ref{Tz}) have a poles of order $\leq 2$ at $P_0$.  Since
$\omega_a(z)$, $G(z,z_\al)$, $G(z_\al,z)$ and
$G_{2\la}(z,z_\al)$ are regular at $z = P_0$, the pole at $P_0$ of the
two last lines of (\ref{Tz}) is of order at most $p$.  \hfill \qed

The result on the rank of the family $(T_z)$ now follows from the fact
that $h^0(\Omega_X^2(pP_0)) = 3g-3 + p$.

Let us prove 2). If we set $p = \infty$ in the result of 1), we see
that the operators $T_z$, $z$ in $U$, commute on all functions of
$\la$ and the $z_i$, which are symmetric in these variables and behave
as sections of $\Omega_X\cL_{-2\la}$, regular outside $P_0$. The
commutator $[T_z,T_z']$ is again a differential-evaluation operator.
But no such operator can vanish on these functions without being zero.
\hfill \qed \medskip

\begin{remark}
  Arguments similar to the proof of Thm.\ \ref{main} imply that the
  $T_z$ defined by (\ref{Tz:hrk}) commute when $k$ is critical.  
\end{remark}

\begin{remark} In the case $n=0$, we find a commuting family of operators
\begin{equation} \label{n=0}
  (T_z f)(\la_1,\ldots,\la_g) = [{1\over 2}(\sum_a \omega_a(z)
  \pa_{\la_a})^2 + \sum_a D_z^{(2\la)}\omega_a(z) \pa_{\la_a} - 2
  \omega_{2\la} (z)] f(\la_1,\ldots,\la_g).
\end{equation} 

If $g = 1$, we have $\omega_a = 2i\pi dz$, $G_\la(z,z') =
{{\theta(-{{\la}\over{2i\pi}} + z -
    z')\theta'(0)}\over{\theta(-{{\la}\over{2i\pi}})\theta(z-z')}}
dz$, $D_z^{(2\la)}\omega_a =
2{{\theta'}\over{\theta}}({{\la}\over{i\pi}}) 2i\pi(dz)^2$,
$\omega_{2\la} = - {{\theta''}\over{\theta}}({{\la}\over{2i\pi}})(dz)^2$,
where $\theta$ is the Jacobi theta-function, so that
\begin{align*} T_z & = [{1\over 2}(2i\pi \pa_\la)^2 + 2
  {{\theta'}\over{\theta}}({{\la}\over{i\pi}}) 2i\pi \pa_\la + 2
  {{\theta''}\over{\theta}}({\la\over{i\pi}})] (dz)^2\\ & = {1\over
    2} (2i\pi \pa_\la + 2 {\theta'\over\theta}({\la\over{i\pi}}))^2
  (dz)^2,
\end{align*}
which is conjugate to ${1\over 2}(2i\pi\pa_{\la})^2$. 

When $g>1$, (\ref{n=0}) is a generating series for one first order and
$3g-3$ second order operators. The linear operator is $\sum_a 2(1-g)
\omega_a(P_0)\pa_{\la_a} + (1-g)a_{00}(\la)$.  From the formula for
the variation of the periods matrix $\delta \tau_{ab} =
\res_{P_0}(\omega_a \omega_b \xi)$ follows that the operator
corresponding to a variation $\delta\tau_{ij}$ has leading term
$\sum_{a,b} \delta\tau_{ab} \pa_{\la_a}\pa_{\la_b}$.
\end{remark}

\begin{remark}
  In the case of the rational curve, we get the commuting family of
  operators defined on symmetric functions $f(z_1,\cdots,z_n)$ by
\begin{equation}
(T(z)f)(z_1,\cdots,z_n) = \sum_{i=1}^n \left( 
\sum_{j\neq i}{1\over{z_j - z_i}} \right)
{{f(z_1,\cdots, z, \cdots,z_n) -f(z_1,\cdots,z_n) }\over{z-z_i}} , 
\end{equation}
where $z$ is at the $i$th position in the r.h.s.
\end{remark}

\begin{remark} \label{Beil-Dr}
  {\it Relation with the Beilinson-Drinfeld operators.} It is not
  possible to interpret the operators $T_z$ directly as
  Beilinson-Drinfeld (BD) operators (\cite{BD}). Indeed, for $g =
  n_{\cK}t[f_\la]w$, with $n_\cK$ in $N(\cK)$, $f_\la$ in $C_\la$ (see
  sect. \ref{tame}) and $w = \pmatrix z^n & 0 \\ 0 &
  z^{-n}\endpmatrix$, the local ring $\hat \cO_{Bun_{\bar G}}([g])$ is
  $H_0(^{g^{-1}}\G^{out}, Ind_{\G_\cO}^{\G}\CC_\chi)^*$, where
  $\CC_\chi$ is the $\G_\cO$-module associated with the character
  $\chi$ of $\G_\cO$, defined by $\chi(K) = -2$ and $\chi(\bar\G
  \otimes \cO) = 0$, and $^{g}x$ denotes the conjugation of $x$ by $g$
  for $x$ in $\G$ and $g$ in $\bar G(\cK)$. This space is isomorphic
  to $CB_\la(W_{-2n|-2n,2n})$, which has no interpretation in terms of
  the $\cF_\la^{(n)}(p)$.

  However, the vector $f[z^{-2n-1}]^{p}v_{-n}$ is cyclic in
  $W_{-2n|-2n,2n}$, which implies that $W_{-2n|-2n,2n}$ is a quotient
  of $W_{p-2n|-2n,2n+2}$. $CB_\la(W_{-2n|-2n,2n})$ may then be viewed
  as a subspace of $CB_\la(W_{p-2n|-2n,2n+2})$, which has a functional
  interpretation when $p\geq 2n$. The BD operators may then be
  expressed as the commuting family of operators $(T_z)$, acting on
  some subspace (defined as the intersection of a family $\wt f[\rho]$
  and some vanishing conditions) of some $\cF_\la^{(n)}(p)$.

  Another connection with the BD operators is the following. The BD
  operators admit lifts to bundles over the moduli space of $\bar
  G$-bundles with parabolic structure at $P_0$.  Such bundles are
  attached to a weight $\La$. The space of local sections of this
  bundle is then $H_0(^{g^{-1}}\G^{out},
  Ind_{\G_{0|0,1}}^{\G}\CC_{\chi_\La})^*$ where $\CC_{\chi_\La}$ is
  the $\G_{0|0,1}$-module defined by $\chi_\La(K) = -2$,
  $\chi_\La(h[1]) = \La$ and $\chi_\La(x[t^i]) = 0$ for $x = f$ and
  $i\ge 0$, and $x = h,e$ and $i>0$. This space is isomorphic to
  $CB_\la(W_{\la -2n | -2n,2n+1})$ which is isomorphic to the
  intersection of kernels of some $\wt f[\rho]$ in some
  $\cF_\la^{(n)}(p)$ if $-2n \leq 1-g$ and $\la>2n$.

  The commuting family of operators $(T_z)$, acting on the
  intersection of kernels of the $\wt f[\rho]$, gets then identified
  with the BD operators. The commuting family $(T_z)$ acting on
  $\cF_\la^{(n)}(p)$ itself gets then identified with the lift of the
  BD operators to some moduli space of $B$-bundles with additional
  structure.  
\end{remark}

\appendix
\section{Proof of Thm.\ \ref{seerose}}

\subsection{Adelization}

For any point $s$ of $X$, denote by $\cK_s$ and $\cO_s$ the local
field and ring at this point. For a finite subset $S$ of $X$, set
$\cK_S = \oplus_{s\in S}\cK_s$ and $\cO_S = \oplus_{s\in S}\cO_s$. Set
also $R_S = H^0(X-S,\cO_X)$; we view $R_S$ as a subring of $\cK_S$.
Define $\G_S$ as the Lie algebra $(\bar\G \otimes \cK_S) \oplus \CC K$,
endowed with the Lie bracket
\begin{equation} \label{Lie}
[x[\eps],y[\eps']] = [x,y][\eps\eps'] + K\langle d\eps ,\eps'\rangle,  
\end{equation} 
with $\langle \omega, \eps\rangle = \sum_{s\in S} \res_s(\omega\eps)$
and $x[\eps] = (x\otimes\eps,0)$. Set $\G_S^{out} = \bar\G \otimes
R_S$; we view $\G_S^{out}$ as a Lie subalgebra of $\G_S$, by the
embedding $x\otimes r \mapsto x[r]$.  For any $s$ in $X$, let
$\G_s$ be the space $(\bar\G \otimes \cK_s) \oplus \CC K$, endowed
with the bracket analogous to (\ref{Lie}), is a Lie subalgebra of
$\G^\AAA$; the associated embedding is denoted by $i_s$.

Let $k$ be a positive integer, $(\La,k)$ be an integrable weight of
$\G$ and $(\rho_{\La,k},L_{\La,k})$ be the associated integrable
module over $\G$.

Define $(\rho_{0,k},L_{0,k})$ as the integrable module over $\G$ with
highest weight $(0,k)$ (the vacuum module of level $k$). Denote by
$v_{top}$ its highest weight vector.  Define $V^S$ as the vector space
$L_{\La,k} \otimes \otimes_{s\in S, s\neq P_0} L_{0,k}$; there is a
map $\rho_S:\G_S \to \End(V^S)$ defined by the condition that the
action of $\G_s$ by $\rho_S \circ i_s$ on $V^S$ is identical to
$\rho_{\La,k}^{(P_0)}$ if $s = P_0$ and to $\rho_{0,k}^{(s)}$ else.

Define $\G^\AAA$ as the space $(\bar\G \otimes \AAA) \oplus \CC K$,
endowed with the Lie bracket analogous to (\ref{Lie}); the map
$x\mapsto (x,0)$ makes $\bar\G \otimes \CC(X)$ a Lie subalgebra of
$\G^\AAA$.  For $x$ in $\bar\G$, $\eps = (\eps_s)_{s\in X}$ in $\AAA$,
we sometimes denote by $x^{(s)}[\eps]$ the element of $\G_s$ equal to
$(x\otimes \eps_s,0)$.

Define $V^{\AAA}$ as the $\G^{\AAA}$-module $\otimes'_{x\in X} V_x$,
with $V_x = L_{0,k}$ for $x\neq P_0$ and $V_{P_0} = L_{\La,k}$.  (Here
$\otimes'$ means that the module is spanned by the products
$\otimes_{x\in X} v_x$ with $v_x$ in $V_x$ equal to the vacuum vector
$v_{top}^{(x)}$ for all but finitely many $x$.) The proof of the
following Lemma is a variant of that of \cite{TUY}, Prop.\ 2.2.3:

\begin{lemma} \label{finkel}
  Let $\psi$ be a $\G^{out}$-invariant linear form on $L_{\La,k}$.
  For any finite subset $S$ of $X$ containing $P_0$, there is a unique
  linear form $\psi_S$ on $V^S$, which is $\G_S^{out}$-invariant and
  such that $\psi_S( \otimes_{x\in S, x\neq P_0} v^{(x)}_{top} \otimes
  v)= \psi(v)$ for any $v$ in $L_{\La,k}$. 

  There is also a unique linear form $\psi^{\AAA}$ on $V^{\AAA}$,
  which is $\bar\G \otimes \CC(X)$-invariant and such that
  $\psi^{\AAA} (\otimes_{x\in X, x\neq P_0} v^{(x)}_{top} \otimes v) =
  \psi(v)$ for any $v$ in $L_{\La,k}$.
\end{lemma}

{\em Proof.} Let us set $\G^{out}_{P_0,x} = H^0(X -
\{P_0,x\},\bar\G)$. Let us denote by $W_{0,k}$ the Weyl module $U\G
\otimes_{U\G^{in}}\CC$, where $\CC$ is the $\G^{in}$-module on which
$\bar\G\otimes \cO$ acts by zero and $K$ acts by $k$. Let us prove
that there is a bijective correspondence between

{\sl (i)} the forms $\psi_{P_0}$ on $L_{\La,k}$, which are
$\G^{out}$-invariant,

{\sl (ii)} the forms $\psi_{P_0,x}$ on $W_{0,k}\otimes L_{\La,k}$, which
are $\G^{out}_{P_0,x}$-invariant

and

{\sl (iii)} the forms $\bar\psi_{P_0,x}$ on $L_{0,k} \otimes L_{\La,k}$, which
are $\G^{out}_{P_0,x}$-invariant, 
the correspondence being such that 
$$\psi_{P_0}(v) =
\psi_{P_0,x}(v_{top}\otimes v) = \bar\psi_{P_0,x}(v_{top}\otimes v).
$$
The proof of the general statement of the Lemma is similar. 

Let us construct a form as in {\sl (ii)} from a form as in {\sl (i)}.
Fix a family of functions $(\rho_i)_{i>0}$ in $H^0(X -
\{P_0,x\},\cO_X)$, such that $\rho_i$ has the expansion $z_x^{-i} +
O(1)$ near $x$, and a basis $(x_\al)_{\al \in A}$ of $\bar\G$. Choose
an order of the index set $A$. By the PBW theorem, a basis of
$W_{0,k}$ is given by the $\prod_\al x_\al^{(x)}[\rho_{i_1(\al)}]
\ldots x_\al^{(x)}[\rho_{i_{n(\al)}(\al)}] v_{top}$, for sequences of
integers $n(\al)$ and of indices $i_1(\al) \leq i_2(\al)
\cdots \leq i_{n(\al)}(\al)$, where the product is performed according
to the order of $A$. Set then
$$ \psi_{P_0,x}(\prod_\al x_\al^{(x)}[\rho_{i_1(\al)}] \ldots
x_\al^{(x)}[\rho_{i_{n(\al)}(\al)}] v_{top} \otimes v) =
\psi_{P_0}(\prod_\al{}^{'} x_\al^{(P_0)}[-\rho_{i_{n(\al)}(\al)}]
\ldots x_\al^{(P_0)}[-\rho_{i_1(\al)}] v).
$$ Here $\prod'$ means that the product over all $\al$'s is taken in
the order inverse to the order of $A$. 
We have then 
$$ \psi_{P_0,x}(\prod_{\al \in A}x_\al^{(P_0,x)}[\rho_{i_1(\al)}] \ldots
x^{(P_0,x)}_\al[\rho_{i_{n(\al)}(\al)}] (v_{top} \otimes v)) = 0 ,
$$ for all $v$ in $V_{\La,k}$, if the product is nonempty. Since the
elements of $W_{0,k} \otimes V_{\La,k}$ are combinations of the
$\prod_{\al \in A}x_\al^{(P_0,x)}[\rho_{i_1(\al)}] \ldots
x^{(P_0,x)}_\al[\rho_{i_{n(\al)}(\al)}] (v_{top} \otimes v)$, it follows
that $\psi_{P_0,x}$ is $\G^{out}_{P_0,x}$-invariant.

Let us now show that any form as in {\sl (ii)} is of the type {\sl
  (iii)}. We follow the argument of \cite{FSV}, based on \cite{Finkelbg}.

For any integer $N \geq 2g$, we can construct an element $\rho_{(N)}$
in $H^0(X-\{P_0,x\},\cO_X)$ with the expansions $\rho_{(N)} = z_x^{-1}
+ O(1)$ near $x$ and $\rho_{(N)} = z_{P_0}^{-N} (\al + O(z_{P_0}))$
near $P_0$, with $\al \neq 0$. For that, it suffices to add to
$\rho_1$ some function of $H^0(X-\{P_0\},\cO_X)$.

Fix $\al^\vee$ in the coroot lattice, such that $\langle \al^\vee,
\theta \rangle \neq 0$.  Let $N$ be an integer $\geq 2g$ and of the form
$1 + d\langle \al^\vee, \theta \rangle$, with $d$ integer.

$L_{0,k}$ is the quotient $W_{0,k}/I$, where $I$ is the submodule of
$W_{0,k}$ generated by $e_{\theta}[z_{x}^{-1}]^{k+1}v_{top}$, where
$e_{\theta}$ is the root vector associated to the maximal root
$\theta$. $I$ is isomorphic to some Verma module. From \cite{Finkelbg}
follows that $e^{(x)}_{\theta}[z_{x}^{-1}]$ is surjective on $I$.  One
may use some some element of the form $\exp(h^{(P_0)}[\eps])$, with
$\eps$ in $z_x\CC[[z_x]]$, to conjugate $e_{\theta}^{(x)}[z_{x}^{-1}]$
to $e_{\theta}^{(x)}[\rho_{(N)}]$.  Therefore,
$e_{\theta}^{(x)}[\rho_{(N)}]$ is also surjective on $I$.

Let us now show that $e_{\theta}^{(P_0)}[\rho_{(N)}]$ is locally
nilpotent on $L_{\La,k}$.  $e_{\theta}^{(P_0)}[\rho_{(N)}]$ is
conjugated by some element of the form $\exp(h^{(P_0)}[\eps])$, with
$\eps$ in $z_{P_0}\CC[[z_{P_0}]]$, to $\al e_{\theta}[z_{P_0}^{-N}]$.
Recall that the affine Weyl group contains a translation element
$w_\omega$ associated to any $\omega$ in the coroot lattice; the
action of $w_\omega$ on the nilpotent loop generators is $w_\omega
\cdot e_{\al}^{(P_0)}[f] = e_\al [(z_{P_0})^{\langle \omega,
  \al\rangle}f]$, for $e_\al$ the root vector associated to any root
$\al$. Moreover, the module $L_{\La,k}$ endowed with the composition
of the action of $\G$ with an affine Weyl group automorphism is again
integrable. It follows that the action of $w\cdot
e_{\theta}[z_{P_0}^{-1}]$, for $w$ any affine Weyl group element, is
locally nilpotent. In particular, for $w = w_{-d\al^{\vee}}$, we find
that $e_{\theta}[z_{P_0}^{-N}]$ is locally nilpotent on $L_{\La,k}$, as
well as $e_{\theta}^{(P_0)}[\rho_{(N)}]$.

These two results imply that $\bar\psi_{P_0,x}$ vanishes on $I \otimes
L_{\La,k}$: indeed, any $v,v'$ in $I$ and $L_{\La,k}$, fix $m$ such
that $(e_{\theta}^{(P_0)}[\rho_{(N)}])^m v'$ vanishes; we may write $v
= (-e_{\theta}^{(x)}[\rho_{(N)}] )^m v''$, with $v''$ in $I$.
$\bar\psi_{P_0,x}(v\otimes v')$ is then equal to $\bar\psi(v'' \otimes
(-e_{\theta}^{(P_0)}[\rho_{(N)}])^m v')$, which is zero.  \hfill \qed
\medskip

\subsection{Formula for the tame symbol} \label{tame}

Denote by $\sigma$ the {\it tame symbol} defined in
$(\AAA^\times)^2$ by
$$ \sigma((f_x)_{x\in X}, (g_x)_{x\in X}) = (-1)^{\sum_{x\in X} v_x(f)
  v_x(g)}\prod_{x\in X} g'(x)^{v_x(f)}f'(x)^{-v_x(g)} ;
$$ we fix a coordinate $z_x$ at each point $x$ of $X$ and set
$f_x = z_x^{v_x(f)}(f'(x) + O(z_x))$.  

Fix a lift $i$ of the universal covering $\wt X \to X$ of $X$,
such that the boundary of $i(X)$ is a union of paths $\wt A_a$, $\wt
B_a$ projecting to a standard system $(A_a),(B_a)$ of $a$- and
$b$-cycles. We will identify the local field and ring at any point $x$
of $X$ with the local field and ring at $i(x)$. For $\la =
(\la_a)$ in $\CC^g$, define $C_\la$ as the set of the adeles of the
meromorphic functions $f: \wt X \to \CC^\times$, such that
$f(\gamma_{A_a}z) = f(z)$ and $f(\gamma_{B_a}z) = e^{-\la_a}f(z)$.

We then have 

\begin{lemma} \label{exist}
  a) For any $\la$ in $\CC^g$, $C_{\la}$ is not empty; moreover, we
  can find elements of $C_{\la}$ without any zero or pole on the $\wt
  A_a$.

  b) For $f$ in $\CC(X)^\times$, without any zero or pole on the
  cycles $A_a$, and $f_{\la}$ in $C_{\la}$, we have
  $$ \si(f,f_\la) = e^{\sum_a n_a(f) \la_a} ,
  $$ with $n_a(f) = {1\over{2i\pi}}\int_{A_a}{df\over f}$.
\end{lemma}

{\em Proof.} Let us prove a). Denote by $\Theta$ the Riemann
theta-function on the Jacobian on $X$, and by $A$ the Abel map. Let
$a$ be any vector of the Jacobian of $X$, then the function $z\mapsto
{{\Theta(A(z)+ a - \la/2i\pi)}\over {\Theta(A(z) +a)}}$ belongs to
$C_{\la}$.  That the zero-poles requirement can be satisfied follows
from a transversality argument.

Let us prove b). Suppose that $f,g$ are nonzero meromorphic functions
on $i(X)$, such that $\sum_{x\in X}\res_x {{df}\over f} =
\sum_{x\in X}\res_x {{dg}\over g} =0$. Then we may introduce cuts on
$\wt X$, connecting the zeroes and the poles of $f$, and choose a
determination of $\ln(f)$ which is single-valued along $\pa i(X)$.
The same can be done for $g$. We have then
$$
\sigma(f,g) = \exp( {1\over{4i \pi}} \int_{\pa i(X)} {df\over f} \ln g 
- {dg\over g} \ln f).   
$$ 
This formula may be proved by deforming $\pa i(X)$ to a set of
contours encircling the cuts of $\ln f$ and $\ln g$.  

Then in the case where $f$ and $g$ belong to $\CC(X)^\times$ and
$C_\la$, we evaluate the integral comparing the contributions of the
paths above $A_a$ and $A_a^{-1}$, and above $B_a$ and $B_a^{-1}$. For
example, in case the zeroes and poles of $f$ and $g$ form disjoint
sets, integration by parts gives
\begin{align*}
  {1\over{4i\pi}}\int_{\pa i(X)} {df\over f} \ln g - {dg\over g} \ln
  f & = {1\over{2i\pi}}\int_{\pa i(X)} {df\over f} \ln g \\ & =
  {1\over{2i\pi}}\sum_a \int_{A_a} {{df}\over f}(\ln g(z) - \ln
  g(\gamma_{B_a}z) ) \\ & = {1\over{2i\pi}}\sum_a \int_a{{df}\over
    f}\la_a,
\end{align*}
which implies b).  \hfill \qed \medskip

\begin{remark}
  Lemma \ref{exist}, b) implies that $\sigma(f,g) = 1$ for any $f,g$
  in $\CC(X)^\times$, which is a well-known fact. One could also
  prove that for any $f$ in $C_\la$ and $f'$ in $C_{\la'}$, without
  any zero or pole on the $\wt A_a$, we have 
\begin{equation} \label{la':la''}
  \sigma(f,f') = e^{\sum_a n_a(f)\la'_a - n_a(f')\la_a}.
\end{equation}
  \hfill \qed \medskip
\end{remark}

\subsection{Construction of $\wt \psi_\la$}

We now follow the classical procedure to construct operators in
$\End(V^{\AAA})$ integrating the Lie algebra action on $V^{\AAA}$.
For $f$ in $\AAA$, $e_i[f]$ and $f_i[f]$ are locally nilpotent on
$V^\AAA$. We set
$$
n^+_i[f] = \exp(e_i[f]), \quad n^-_i[f] = \exp(f_i[f])
$$ for $f$ in $\AAA$. Set also, for $\rho$ in $\AAA^\times$,
$w_i[\rho] = n^+_i[\rho]n^-_i[-\rho^{-1}]n^+_i[\rho]$, and
$$
t_i[\rho] = w_i[\rho]w_i[1]^{-1}. 
$$
We have then 
\begin{equation} \label{t:t}
  t_{i}[\rho\rho'] = \sigma(\rho,\rho')^{-k(h_i|h_i)/2}
  t_{i}[\rho]t_{i}[\rho']
\end{equation}
for $i$ simple, and 
\begin{equation} \label{ti:tj} 
  t_{i}[\rho]t_{j}[\rho']t_{i}[\rho]^{-1} t_{j}[\rho']^{-1} =
  \sigma(\rho,\rho')^{k(h_i|h_j)},
\end{equation} 
for any indices $i,j$ (observe that $(h_i|h_j)$ is always integer and
$(h_i|h_i)$ always even).

The first identity is a consequence of \cite{Garland}, Thm.\ 12.24, and
the second is a consequence of this identity and \cite{Steinbg}, 7.3)
e) (see also \cite{Moore}, Lemma 8.2, formula (3)).

\begin{prop:def}
  Let us fix $\la^{(1)},\ldots, \la^{(r)}$ in $\CC^g$. For
  $f_{\la^{(i)}}$ in $C_{\la^i}$, such that the $f_{\la^{(i)}}$ have
  no zero or pole on the $A_a$, and $v$ in $V_{\La,k}$, the quantity
  \begin{equation} \label{auble}
    e^{ \sum_i {{k(h_i|h_i)}\over {2} } \sum_a \la^{(i)}_a
      n_a(f_{\la^{(i)}}) + \sum_{i< j} k(h_i|h_j) \sum_a \la^i_a
      n_a(f_{\la^{(j)}}) }\langle \psi^{\AAA} , t_1[f_{\la^{(1)}}] \cdots
    t_r[f_{\la^{(r)}}] \left( v \otimes \otimes_{x\neq P_0} v_{top}^{(x)}
    \right) \rangle
  \end{equation} 
  is independent of the choice of the $f_{\la^{(i)}}$. We will set
  $\la = (\la^{(1)},\ldots,\la^{(r)})$ and
  \begin{align*} & \langle \wt\psi_\la , v \rangle  = 
    e^{ \sum_i {{k(h_i|h_i)}\over {2} } \sum_a \la^{(i)}_a
      n_a(f_{\la^{(i)}}) + \sum_{i< j} k(h_i|h_j) \sum_a \la^{(i)}_a
      n_a(f_{\la^{(j)}}) } \cdot \\ & \cdot \langle \psi^{\AAA} ,
    t_1[f_{\la^{(1)}}] \cdots t_r[f_{\la^{(r)}}] \left( v \otimes
      \otimes_{x\neq P_0} v_{top}^{(x)} \right) \rangle
\end{align*} for any such $f_{\la^{(i)}}$. 
\end{prop:def}

{\em Proof.} Let $f'_{\la^{(i)}}$ be other elements of $C_{\la^{(i)}}$,
satisfying the same zero-poles condition as $f_{\la^{(i)}}$. Then
$f'_{\la^{(i)}} = f_if_{\la^{(i)}}$, with $f_i$ in $\CC(X)^\times$, without
zero or pole on the $A_a$.  We have 
\begin{align*}
  & e^{ \sum_i {{k(h_i|h_i)}\over {2} } \sum_a \la^{(i)}_a
    n_a(f'_{\la^{(i)}}) + \sum_{i< j} k(h_i|h_j) \sum_a \la^{(i)}_a
    n_a(f'_{\la^{(j)}}) }\langle \psi^{\AAA} , t_1[f'_{\la^{(1)}}]
  \cdots t_r[f'_{\la^{(r)}}] \left( v \otimes \otimes_{x\neq P_0}
    v_{top}^{(x)} \right) \rangle \\ & = e^{ \sum_i {{k(h_i|h_i)}\over
      {2} } \sum_a \la^i_a n_a(f_i) + \sum_{i< j} k(h_i|h_j) \sum_a
    \la^{(i)}_a n_a(f_j) } \cdot \\ & \cdot e^{ \sum_i
    {{k(h_i|h_i)}\over {2} } \sum_a \la^{(i)}_a n_a(f_{\la^{(i)}}) +
    \sum_{i< j} k(h_i|h_j) \sum_a \la^{(i)}_a n_a(f_{\la^{(j)}}) }
  \cdot \\ & \cdot \langle \psi^{\AAA} , t_1[f_1f_{\la^{(1)}}] \cdots
  t_r[f_r f_{\la^{(r)}}] \left( v \otimes \otimes_{x\neq P_0}
    v_{top}^{(x)} \right) \rangle
\end{align*}
The identities (\ref{t:t}) and (\ref{ti:tj}) imply that this is
equal to
\begin{align} \label{tran}
  & \nonumber e^{ \sum_i {{k(h_i|h_i)}\over {2} } \sum_a \la^{(i)}_a
    n_a(f_i) + \sum_{i< j} k(h_i|h_j) \sum_a \la^i_a n_a(f_j) } \\ &
  \nonumber e^{ \sum_i {{k(h_i|h_i)}\over {2} } \sum_a \la^i_a
    n_a(f_{\la^{(i)}}) + \sum_{i< j} k(h_i|h_j) \sum_a \la^i_a
    n_a(f_{\la^{(j)}}) } \cdot \\ & \nonumber \cdot \prod_i
  \si(f_i,f_{\la^{(i)}})^{- {{k(h_i|h_i)}\over 2}} \prod_{i<j}
  \si(f_j,f_{\la^{(i)}})^{- k(h_i|h_i)} \\ & \langle \psi^{\AAA} ,
  t_1[f_1]\cdots t_r[f_r] t_1[f_{\la^{(1)}}] \cdots t_r[f_{\la^{(r)}}]
  \left( v \otimes \otimes_{x\neq P_0} v_{top}^{(x)} \right) \rangle
\end{align}

Now, as the $t_i[f_i]$ are products of exponentials of elements of the
$\bar\G\otimes \CC(X)$ and $\psi^\AAA$ is $\bar\G \otimes
\CC(X)$-invariant, we have $\langle \wt\psi_\la, \prod_{i=1}^r t_i[f_i]
v'\rangle = \langle \wt\psi_\la, v'\rangle$ for any $v'$ in $V^\AAA$.
Applying Lemma \ref{exist}, b), we find that (\ref{tran}) is equal to
$$ e^{  \sum_i {{k(h_i|h_i)}\over {2} } \sum_a \la^{(i)}_a
  n_a(f_{\la^{(i)}}) + \sum_{i< j} k(h_i|h_j) \sum_a \la^{(i)}_a
  n_a(f_{\la^j}) }\langle \psi^{\AAA} , t_1[f_{\la^{(1)}}] \cdots
t_r[f_{\la^{(r)}}] \left( v \otimes \otimes_{x\neq P_0} v_{top}^{(x)}
\right) \rangle .
$$ \hfill \qed 

\begin{remark} \label{megido}
  In view of (\ref{ti:tj}) and (\ref{la':la''}), it is clear that
  (\ref{auble}) is independent of the chosen ordering of simple
  coroots.
\end{remark}

Let us now give an expression of $\wt\psi_\la$ in terms of extremal
vectors.

\begin{lemma} \label{5.6}
  Define the vectors $v_{i;[n]}$ in $L_{0,k}$ by the formulas
  $v_{i;[0]} = v_{top}$, $v_{i;[n+1]} =
  {{(-1)^k}\over{k!}}f_i[z^{-2n-1}]^k v_{i;[n]}$ and $v_{i;[-n-1]} =
  {{1}\over{k!}}e_i[z^{-2n-1}]^k v_{i;[-n]}$ for $n\geq 0$. Then we
  have $v_{i;[n]} = t_i[z^n]v_{top}$.
\end{lemma}

{\em Proof.} It is enough to prove this statement for the case $\bar\G
= \SL_2$. The formulas for $v_{i;[1]}$ and $v_{i;[-1]}$ are derived by
direct expansions. The other formulas are obtained by applying the
affine Weyl group translation associated with the coroot $h_i$ (which
preserves $t_i[z]$). \hfill \qed

We have then 

\begin{prop} \label{extremal}
  Assume that the sets $S_i$ of zeroes and poles of the $f_{\la^{(i)}}$ are
  distinct.  Then we have for $v$ in $V_{\La,k}$,
\begin{align*}
  & \langle \wt\psi_\la , v \rangle = e^{ \sum_i {{k(h_i|h_i)}\over
      {2} } \sum_a \la^{(i)}_a n_a(f_{\la^{(i)}}) + \sum_{i< j}
    k(h_i|h_j) \sum_a \la^{(i)}_a n_a(f_{\la^{(j)}}) } \\ &
  \prod_{i=1}^r \prod_{s\in S_i} (f'_{\la^{(i)}}(s)^{-k}) \langle
  \psi_{\{P_0\}\cup (\cup_i S_i)} , \otimes_{i=1}^r (\otimes_{s\in
    S_i} v^{(s)}_{i;[v_{s}(f_{\la^{(i)}})]} ) \otimes \prod_{i=1}^r
  t_i^{(P_0)}[f_{\la^{(i)}}] v \rangle ,
\end{align*}
where we set $f_{\la^{(i)}}(z) = f'_{\la^{(i)}}(s) z_s + o(z_s)$ for
$s$ in $S_i$. Recall that $\psi_{\{P_0\}\cup (\cup_i S_i)}$ denotes
the prolongation of $\psi$ to the product of $L_{\La,k}$ and vacuum
modules at the points of $S_i$.
\end{prop}

\subsection{Proof of Thm.\ \ref{seerose}}

To prove Thm.\ \ref{seerose}, 1), we first prove 
\begin{lemma} \label{mahaz} 
  For any $v$ in $L_{\La,k}$, the function $\la\mapsto \langle
  \wt\psi_\la , v \rangle$ depends analytically on $\la$ and satisfies 
  $$ \pa_{\la^{(i)}_a} \langle \wt \psi_\la , v\rangle = \langle \wt
  \psi_\la , h_i[r_a]v\rangle,
  $$ $a = 1, \ldots,g, i = 1, \ldots, r$.
\end{lemma}

{\em Proof of Lemma.}
Let us prove this first in the case $\bar \G = \SL_2$.  In that case, we
work in a neighborhood of some point $\la_0$ of $J^0(X)$.  Let
$P_i(\la)$ be points on $X$ ($i = 1, \ldots,g$) such that $f_\la$ has
simple zeroes at the $P_i(\la)$ and a pole of order $g$ at $P_0$. Let
$z_{P_i(\la_0)}$ a coordinate at $P_i(\la_0)$; we will again denote by
$P_i(\la)$ the coordinate of the point $P_i(\la)$ in the coordinate
system. We will assume that the local coordinate at $P_i(\la)$ is
$z_{P_i(\la)} = z_{P_i(\la_0)} - P_i(\la)$.

Let for $P$ in $X$, $\rho_{P}$ be a meromorphic function on $X$, with
only poles at $P_0$ and at $P$, with the expansion $\rho = z_P^{-1} +
O(1)$. We assume that the expansions at $P_0$ the functions
$\rho_P$ depend smoothly on $P$, for $P$ near any of the $P_i(\la_0)$.
We set also $f_\la(z) = f'_\la(P_i(\la))z_{P_i(\la)} + {1\over
  2}f''_\la(P_i(\la))z_{P_i(\la)}^2 + \cdots$.

Then Prop.\ \ref{extremal} implies that 
$$ \langle \wt\psi_\la , v \rangle = e^{k \sum_a \la_a n_a(f_\la)}
\prod_{i} f'_\la(P_i(\la))^{-k}
\langle \psi_{ \{P_0,P_i(\la)\} }, \otimes_{i=1}^g
{{f[-z_{P_i(\la)}^{-1}]^k}\over{k!}}v^{(P_i(\la))}_{top} \otimes
t^{(P_0)}[f_\la] v \rangle .
$$ As we have seen, $f[-z_{P_i(\la)}^{-1}]^k v^{(P_i(\la))}_{top}$ is
equal to $f^{(P_i(\la))}[-\rho_{P_i(\la)}]^k
v^{(P_i(\la))}_{top}$.  By the coinvariance of $\psi$, and the fact
that $v_{[1]}$ is annihilated by the $f[\phi]$, $\phi$ in $\cO$, the
r.h.s. of this equation is equal to
$$ {1\over{(k!)^g}}e^{k \sum_a \la_a n_a(f_\la)} \prod_{i}
f'_\la(P_i(\la))^{-k}\langle \psi, \prod_{i=1}^g \left(
  f[\rho_{P_i(\la)}]^k \right) t[f_\la^{(P_0)}] v \rangle .
$$ This formula shows that $\langle \wt\psi_\la , v\rangle$ depends
smoothly on $\la$. Let us compute its differential. Let $\delta\la$ be
a variation of $\la$.  A computation of adjoint actions shows that
$$ \delta t[f_\la^{(P_0)}] = (h[{{\delta
    f_\la^{(P_0)}}\over{f_\la^{(P_0)}}}] + k \langle df_\la^{(P_0)},
{{\delta f^{(P_0)}_\la}\over{(f^{(P_0)}_{\la})^2}}\rangle) t[f_\la^{(P_0)}], 
$$ so that
\begin{align} \label{tonio}
  & \nonumber \delta \langle \wt\psi_\la,v \rangle = k(\sum_a
  n_a(f_\la)\delta \la_a) \langle \wt\psi_\la,v \rangle \\ & \nonumber +
  {1\over{(k!)^g}}e^{k \sum_a \la_a n_a(f_\la)} \sum_{i=1}^g \langle
  \psi, \prod_{j\neq i} \left( f[\rho_{P_j(\la)}]^k \right) k
  f[\rho_{P_i(\la)}]^{k-1} f[\delta\rho_{P_i(\la)}] t[f_\la^{(P_0)}] v
  \rangle \\ & \nonumber \prod_{i} f'_\la(P_i(\la))^{-k}\\ & \nonumber
  + {1\over{(k!)^g}}e^{k \sum_a \la_a n_a(f_\la)} \langle \psi,
  \prod_{i=1}^g \left( f[\rho_{P_i(\la)}]^k \right) (h[{{\delta
      f_\la^{(P_0)}}\over{f_\la^{(P_0)}}}] + k \langle df_\la^{(P_0)},
  { {\delta f^{(P_0)}_\la}\over{(f^{(P_0)}_{\la})^2}}\rangle)
  t[f_\la^{(P_0)}] v \rangle \\ & \nonumber \prod_{i}
  f'_\la(P_i(\la))^{-k } \\ & \nonumber + (-k \sum_{i=1}^g {{\delta
      f'_\la(P_i(\la))}\over{f'_\la(P_i(\la))}}) \langle \wt\psi_\la,
  v\rangle ,
\end{align}
which can be rewritten (using coinvariance) as
\begin{align*}
  & \delta \langle \wt\psi_\la, v\rangle = (k \sum_a \delta\la_a
  n_a(f_\la) -k \sum_i {{\delta
      f'_{\la}(P_i(\la))}\over{f'_\la(P_i(\la))}} + k \langle
  {{df_\la}\over{f_\la}}, {{\delta f_\la}\over{f_\la}}\rangle_{P_0} )
  \langle \wt\psi_\la , v \rangle \\ & + \prod_{i} f'_\la(P_i(\la))^{-k}
  \\ & \langle \psi_{\{P_0,P_i(\la)\}}, \sum_i \otimes_{j\neq i}
  v_{[1]}^{(j)} \otimes {{(-1)^k}\over{k!}} f[\delta
  P_i(\la)z_{P_i(\la)}^{-2}]f[z^{-1}_{P_i)(\la)}]^{k-1}v_{top}^{(i)}
  \otimes t^{(P_0)}[f_\la] v \rangle \\ & + \prod_{i}
  f'_\la(P_i(\la))^{-k} \langle \psi_{\{P_0,P_i(\la)\}}, \sum_i
  h^{(P_0)}[{{\delta f_\la}\over{f_\la}}] (\otimes_{i=1}^g
  v_{[1]}^{(i)} \otimes t^{(P_0)}[f_{\la}] v) \rangle .
\end{align*}

The penultimate term is rewritten as
$$ \prod_{i} f'_\la(P_i(\la))^{-k} \langle \psi_{\{P_0,P_i(\la)\}},
 - \sum_i \delta P_i(\la) h^{(i)}[z_{P_i(\la)}^{-1}] (\otimes_{i=1}^g
v_{[1]}^{(i)}) \otimes t^{(P_0)}[f_\la] v \rangle ,
$$
using the identity in $L_{0,k}$
$$ h[z^{-1}]v_{[1]} = {{(-1)^{k-1}}\over{(k-1)!}}
f[z^{-2}]f[z^{-1}]^{k-1}v_{top},
$$
which follows from 
\begin{equation} \label{vx}
h[z^{-1}]f[z^{-1}]^k v_{top} = 
\left( e[z]f[z^{-1}] - f[z^{-2}]e[z]\right) f[z^{-1}]^k v_{top} = 
-kf[z^{-2}]f[z^{-1}]^{k-1} v_{top}, 
\end{equation} 
because $f[z^{-2}]f[z^{-1}]^k v_{top} = 0$, which is a consequence of
the integrability conditions.

On the other hand, we have $t[f_{\la}^{(P_0)}] h[r_a]
t[f_{\la}^{(P_0)}]^{-1} = h[r_a] + 2k \langle
{{df^{(P_0)}_\la}\over{f^{(P_0)}_\la}}, r_a \rangle$ so that
$\sum_a \delta \la_a \langle \wt\psi_\la , h[r_a] v \rangle$ is equal to
\begin{align*} 
  & {1\over{(k!)^g}}e^{\sum_a \la_a n_a(f_\la)} \prod_{i}
  f'_\la(P_i(\la))^{-k} \\ & \sum_a \delta \la_a \langle \psi ,
  \prod_{i=1}^g \left( f[\rho_{P_i(\la)}]^k \right) (h[r_a] + 2k
  \langle {{df^{(P_0)}_\la}\over{f^{(P_0)}_\la}}, r_a \rangle)
  t[f_\la^{(P_0)}] v \rangle.
\end{align*}

Therefore, we have 
\begin{align*} 
  & \delta\langle \wt\psi_\la , v\rangle - \sum_a \delta\la_a \langle
  \wt\psi_\la , h[r_a]v\rangle \\ & = [k \sum_a \delta \la_a
  n_a(f_\la) + k\langle {{df_\la}\over{f_\la}} , {{\delta
      f_\la}\over{f_\la}}\rangle_{P_0} -2k \sum_a \delta \la_a \langle
  {{df_\la}\over{f_\la}}, r_a \rangle_{P_0} - k \sum_{i} {{\delta
      f'_\la(P_i(\la))}\over{f'_\la(P_i(\la))}}] \\ & \langle
  \wt\psi_\la , v\rangle \\ & + \prod_{i} f'_\la(P_i(\la))^{-k} \\ &
  \langle \psi_{\{P_0,P_i(\la)\}}, (h^{(P_0)}[{{\delta
      f_\la}\over{f_\la}}] - \sum \delta
  P_i(\la)h^{(i)}[z_{P_i(\la)}^{-1}] - \sum_a \delta\la_a
  h^{(P_0)}[r_a] ) (\otimes_i v_{[1]}) \otimes t[f_\la^{(P_0)}]v
  \rangle .
\end{align*}
On the other hand, $\varrho = {{\delta
    f_\la^{(P_0)}}\over{f_\la^{(P_0)}}} - \sum_a \delta\la_a r_a$ is
single-valued on $X$ and has simple poles at the $P_i(\la)$. 
Therefore, 
$$ \langle \psi_{\{P_0,P_i(\la)\}}, (h^{(P_0)}[\varrho] +
\sum_{i}h^{(i)}[\varrho]) ( (\otimes_i v^{(i)}_{[1]}) \otimes
t^{(P_0)}v)\rangle
$$ is zero, so that $\delta\langle \wt\psi_\la , v\rangle - \sum_a
\delta\la_a \langle \wt\psi_\la , h[r_a]v\rangle$ is proportional to 
\begin{align} \label{sch}
  & - \sum_{i=1}^g {{\delta f'_\la(P_i(\la))}\over{f'_\la(P_i(\la))}}
  + \sum_a \delta\la_a n_a(f_\la) + \langle {{df_\la}\over{f_\la}},
  {{\delta f_\la}\over{f_\la}}\rangle_{P_0} \\ & \nonumber - 2 \sum_a
  \delta \la_a\langle {{df_\la}\over{f_\la}},r_a \rangle_{P_0} + 2
  \sum_{i}\left[ ({{\delta f_\la}\over{f_\la}})^{reg}(P_i(\la)) -
    \sum_a \delta\la_a r_a(P_i(\la)) \right] ,
\end{align} 
where we set $({{\delta f_\la}\over{f_\la}})(z) =
\al_{\la,i}z_{P_i(\la)}^{-1} +  ({{\delta
    f_\la}\over{f_\la}})^{reg}(P_i(\la)) + O(z_{P_i(\la)})$. The
vanishing of (\ref{sch}) then follows from the identities
$$ \langle {{df_\la}\over{f_\la}},{{\delta f_\la}\over{f_\la}}
\rangle_{P_0} = -\sum_i \langle {{df_\la}\over{f_\la}},{{\delta
    f_\la}\over{f_\la}} \rangle_{P_i} + \sum_a n_a(f_\la) \delta\la_a , 
$$
$$ \langle {{df_\la}\over{f_\la}}, r_a\rangle_{P_0} = - \sum_{i}
r_a(P_i(\la)) + n_a(f_\la)
$$
and
$$ - {{\delta f'_\la(P_i(\la))}\over{f'_\la(P_i(\la))}} - \langle
{{df_\la}\over{f_\la}} \rangle_{P_i} + 2 ({{\delta
    f_\la}\over{f_\la}})^{reg}(P_i) = 0 ;  
$$
the latter identity follows from the expansions 
$$ {{df_\la}\over{f_\la}} = {{dz}\over{z-P_i(\la)}} + {1\over
  2}{{f''_\la}\over{f'_\la}}(P_i(\la)) dz + O(z - P_i(\la))dz,
$$
$$ {{\delta f_\la}\over{f_\la}} = - {{\delta
    P_i(\la)}\over{z-P_i(\la)}} + [{{\delta
    f'_\la(P_i(\la))}\over{f'_\la(P_i(\la))}} - {1\over 2} \sum_a
{{f''_\la}\over{f'_\la}}(P_i(\la))\delta P_i(\la)] + O(z - P_i(\la)), 
$$
$$ ({{\delta f_\la}\over{f_\la}})^{reg}(P_i(\la)) = {{\delta
    f'_\la(P_i(\la))}\over{f'_\la(P_i(\la))}} - {1\over 2} \sum_i
{{f''_\la}\over{f'_\la}}(P_i(\la)) \delta P_i(\la).
$$ This ends the proof of Lemma \ref{mahaz} in the case $\bar \G =
\SL_2$. In the case of general $\bar\G$, this result allows to compute
$\pa_{\la_a^{(1)}} \langle \wt\psi_{\la} , v \rangle$; the additional
prefactors of the expression of $\langle \wt\psi_\la, v\rangle$ allow
to transfer the $h_1[r_a]$ in front of $v$. Using Rem. \ref{megido},
we can treat the case of any simple coroot in the same way.
\hfill \qed \medskip 

Let us now show why Lemma \ref{mahaz} implies Thm.\ \ref{seerose}, 1).
The differential equation of Lemma \ref{mahaz} and the equality
$\wt\psi_0 = \psi$ imply that the formal expansion of $\langle
\wt\psi_\la, v\rangle$ for $\la$ near $0$ is equal to $\langle
\psi_\la, v\rangle$. This implies Thm.\ \ref{seerose}, 1).

Thm.\ \ref{seerose}, 2) follows from the equality $\psi_\la =
\wt\psi_\la$ and the fact that for any $f_{\la^{(i)}}$ in
$C_{\la^{(i)}}$, we have $$ \Ad( t_1[f_{\la^{(1)}}]\cdots
t_r[f_{\la^{(r)}}]) (\G^{out}_\la) = \G^{out}.$$

Finally, Thm.\ \ref{seerose}, 3) follows from the equality $\wt\psi_\la
= \psi_\la$ and the fact that if $f_\la$ belongs to $C_\la$, $f_\la
e^{\zeta_a}$ belongs to $C_{\la + \Omega_a}$.  
This ends the proof of Thm.\ \ref{seerose}.

\begin{remark}
  Equation (\ref{vx}) is translated through the states-fields
  correspondence into the identity
  $$ {d\over{dz}}(f(z)^k) = - :h(z) f(z)^k:,
  $$ which is valid in level $k$ modules (see \cite{Lep-Primc}), and
  means that $f(z)^k$ is a vertex operator. The connection between
  this vertex algebra and the Abel-Jacobi map was noticed in
  \cite{FS}.
\end{remark}



\end{document}